\documentclass[12pt,oneside]{amsart}
\usepackage{amsfonts,amssymb,amscd,amsthm,amsmath}
\usepackage{graphpap}

 

\newtheorem{theorem}[subsection]{Theorem}
\newtheorem*{theorem*}{Theorem}
\newtheorem{lemma}[subsection]{Lemma}
\newtheorem{proposition}[subsection]{Proposition}
\newtheorem{corollary}[subsection]{Corollary}

\theoremstyle{definition}
\newtheorem{definition}[subsection]{Definition}
\newtheorem*{problem}{Problem}
\newtheorem{example}[subsection]{Example}
\theoremstyle{remark}
\newtheorem{remark}[subsection]{Remark}

\newcommand{\mt}[1]{\operatorname{#1}}
\newcommand{\EEE}{{\mathbb E}}
\newcommand{\DDD}{{\mathbb D}}
\newcommand{\AAA}{{\mathbb A}}
\newcommand{\QQ}{{\mathbb Q}}
\newcommand{\ZZ}{{\mathbb Z}}
\newcommand{\CC}{{\mathbb C}}
\newcommand{\OO}{{\mathcal O}}

\newcommand{\PP}{{\mathbb P}}

\newcommand{\TT}{{\mathcal T}}
\newcommand{\NN}{{\mathbb N}}
\newcommand{\FFF}{{\mathbb F}}
\newcommand{\an}{{\cong_{\mt {an}}}}
\newcommand{\SM}{{\Phi_{\bf sm}}}
\newcommand{\M}{{\Phi_{\bf m}}}
\newcommand{\Supp}{\mt{Supp}}
\newcommand{\Sing}{\mt{Sing}}
\newcommand{\Diff}{\mt{Diff}}

\newcommand{\Exc}{\mt{Exc}}

\newcommand{\down}[1]{\llcorner #1 \lrcorner}

\title{Classification of exceptional log Del Pezzo surfaces with {\large $\delta=1$}}
\author{S.~A.~Kudryavtsev}

\date{}
\address{Department of Algebra, Faculty of Mathematics,
Moscow State Lomonosov University, 117234 Moscow, Russia}

\email{kudryav@mech.math.msu.su}

\begin{document}
\begin{abstract} The exceptional log Del Pezzo surfaces with
$\delta=1$ are classified.
\end{abstract}
\maketitle

\section*{\bf {Introduction}}

Now the main problem in log Minimal Model Program is the study of
extremal contractions, singularities. In solving this problem
the first step applied in the three-dimensional Fano manifold classification
is to find
a "good"\ divisor in the multiple anticanonical linear system.
Recently V.V.~Shokurov have suggested the inductive method to
construct a "good"\ divisor $D\in |-nK_X|$ in the paper \cite{Sh1}.
Such divisor is called $n$-complement.
The essence of method is the following one: the construction of
$n$-complement for $m$-dimensional extremal contraction $X\to Z$ of local type
(that is $\dim Z>0$) or singularity
(that is $Z=X$) follows from the construction of $n$-complement for
$(m-1)$-dimensional  log variety. Let us demonstrate it in the case of
singularities.
\par
Let $(X\ni P)$ be a klt singularity.
Then there exists a plt blow-up
$f\colon (Y,S)\to (X\ni P)$ \cite[theorem 1.5]{Kud2}, i.e.
$S=\Exc f$ is an irreducible, $f$\ anti-ample divisor and $(Y,S)$ is a plt pair.
Then $n$-complement of $K_S+\Diff_S(0)$ is extended to
$n$-complement of $K_Y+S$ and its image is
$n$-complement of $K_X$.
The converse statement is also true, i.e.
$n$-complement of $K_X$ induces $n$-complement of $K_Y+S$ and
$K_S+\Diff_S(0)$. It should be noticed that the coefficients of different
$\Diff_S(0)$ are standard, i.e. they have form $\frac{k_i-1}{k_i}$, where
$k_i\in \ZZ_{>0}\cup\{\infty\}$.
See \cite{Kud1} in the case of hypersurface singularities.
Also see examples \ref{ex1}--\ref{ex3}.
\par
Therefore there appears an important problem when studying the three-dimensional
singularities and contractions.
\begin{problem} Describe the class of all log Del Pezzo surfaces,
$\PP^1$ and elliptic fibrations which can be exceptional divisors of
some plt blow-ups of three-dimensional log canonical singularities.
\end{problem}

When finding complement the next important concept is an exceptionality.
The exceptionality also remains valid after an inductive transfer from
an extremal contraction $X\to Z$ of local type to a log variety $(S,\Diff_S(0))$.
The importance of separation of extremal contractions, log varieties on
exceptional and non-exceptional ones follows from following two properties.

\begin{enumerate}
\item If a variety or extremal contraction is nonexceptional
then the linear system $|-nK_X|$ must have a "good"\ member for small $n$.
For example, we can take $n\in \{1,2\}$ for two-dimensional
log canonical singularities
\cite[example 5.2]{Sh1} and $n\in \{1,2,3,4,6\}$ for three-dimensional log canonical
singularities \cite[theorem 7.1]{Sh2}.
\item The exceptional singularities are "bounded"\ and can be classified in details.
\end{enumerate}

\par
The boundedness relates to some invariants, for instance to
minimal log discrepancies (see \cite[\S 7]{Sh2}). Hypothetically,
the boundedness of exceptional singularities means the following:
$(S,\Diff_S(0))$ belongs to the finite number of algebraic families.
For three-dimensional exceptional log canonical singularities it was checked in
\cite{Kud1}. The other results about a boundedness and references about it can be
found in \cite{PrLect} and \cite{PrSh}.
\par
Therefore to classify the three-dimensional extremal contractions and
singularities we have to know the classification of the following exceptional log
surfaces $(S,D)$, where $D$ is a boundary with standard coefficients,
$-(K_S+D)$ is an ample divisor and $K_S+D^+$ is a klt divisor for any
complement $D^+$ of $K_S+D$.
\par
Let us consider the following invariant
\begin{eqnarray*}
\delta(S,D)=\#\Big\{E\mid E\  \text{is a divisor with discrepancy}\
a(E,D)\le -\frac67   \Big\}.
\end{eqnarray*}
\par
By theorem 5.1 \cite{Sh2} we have $\delta(S,D) \le 2$.
The classification of surfaces with $\delta(S,D)=2$ is obtained in \cite[\S 5]{Sh2}.
\par
The main result of this paper is to classify the surfaces in the case
$\delta(S,D)=1$.
These results imply the classification of log Enriques surfaces
with $\delta(S,D)=1,2$ (see \cite{Kud4}, \cite{Kud5}).
\par
Let us remark that the surfaces with $\delta(S,D)=1$ were studied earlier.
The case of "elliptic curve"\ was developed in the preprint \cite{Abe1}.
The remaining cases were studied in \cite{Abe2}. The answer obtained in this
dissertation wasn't justified and wasn't correct.
\par
To study the exceptional log Del Pezzo surfaces another approach was given in
\cite{KeM}.
\par
This paper is organized in the following way. In chapter
\S 1 the main definitions and preliminary results are collected.
In chapter \S 2 the classification theorem is formulated and its corollaries are
proved. In chapter \S 3 the basic definitions and constructions are introduced to
prove the classification theorem.
The classification process is completed in chapters
\S 4, \S 5, \S 6 and \S 7.
\par
I am grateful to Professor  Yu.G.~Prokhorov for valuable remarks.
\par
The research was partially
supported by a grant 02-01-00441 from the Russian Foundation of Basic Research
and a grant INTAS-OPEN 2000\#269.

\section{\bf {Preliminary facts and results}}
All varieties are algebraic and are assumed to be defined over
$\CC$, the complex number field.
The main definitions, terminology and notations used in the paper are
given in \cite{Koetal}, \cite{PrLect}.

\subsection{} {\it List of notations.}
\par
The zero and minimal sections of ruled surface
$\FFF_n=\PP(\OO_{\PP^1}\oplus\OO_{\PP^1}(n))$ are denoted by
$E_0$ and $E_{\infty}$ respectively. Its fiber is denoted by $f$.
\par
The coordinates of weighted projective space $\PP=\PP(a_1,\ldots,a_n)$ are denoted
by $x_1$,\ldots,$x_n$.
The general hypersurface of degree $d$ in $\PP$ is denoted by $X_d$.
In many cases it is enough to require the irreducibility and reducibility of $X_d$.
\par
A log Del Pezzo surface with $\rho(S)=1$ is denoted by
$S(\AAA_1+\AAA_5)$. In the brackets its singularities are written.
In our case the surface has
$\AAA_1$ and $\AAA_5$ singularities.
\par
Put $\SM=\{1-1/m\mid m\in\ZZ_{>0}\cup\{\infty\}\}$ and
$\M=\SM\cup [6/7,1]$. The coefficient $d$ is called {\it standard}, if $d\in\SM$.
\par
The symbol $\an$ means an analytical isomorphism.
\par
The intersection index of $D_1$ and $D_2$
at the point $P$ is denoted by $(D_1\cdot D_2)_P$.
\par
The arithmetical genus of $D$ is denoted by $p_a(D)$.

\begin{definition}
Let $(X/Z,D)$ be a log pair, where $D$ is a subboundary. Then a
\textit{$\QQ$-complement} of $K_X+D$ is a
log divisor $K_X+D'$ such that $D'\ge D$, $K_X+D'$ is lc and
$n(K_X+D')\sim 0$ for some $n\in\NN$.
\end{definition}

\begin{definition}\label{defcompl}
Let $X$ be a normal variety and let $D=S+B$ be a subboundary on
$X$ such that $B$ and $S$ have no common components, $S$ is an
effective integral divisor and $\down{B}\le 0$. Then we say that
$K_X+D$ is \textit{$n$-complementary} if there is a $\QQ$-divisor
$D^+$ such that
\begin{enumerate}
\item
$n(K_X+D^+)\sim 0$ (in particular, $nD^+$ is an integral divisor);
\item
$K_X+D^+$ is lc;
\item
$nD^+\ge nS+\down{(n+1)B}$.
\end{enumerate}
In this situation the
\textit{$n$-complement} of
$K_X+D$ is $K_X+D^+$. The divisor $D^+$ is called an
\textit{$n$-complement} too.
\end{definition}

\begin{theorem}\label{contr} \cite[theorem 3.1]{Sh2}
Let $(X/Z\ni P,D)$ be a log surface of local type
$($i.e. $Z$ is not point$)$, where $f\colon X\to Z\ni P$ is a contraction.
Assume that $-(K_X+D)$ is $f$-nef and $K_X+D$ is lc.
Then there exists an $1$-,$2$-,$3$-,$4$- or $6$-complement of $K_X+D$ near
$f^{-1}(P)$.
\end{theorem}

\begin{definition}
Let $(X/Z\ni P,D)$ be a contraction of varieties, where $D$ is a boundary.
\begin{enumerate}
\item Assume that $Z$ is not a point (local case). Then $(X/Z\ni P,D)$
is said to be \textit{exceptional}
over $P$ if for any $\QQ$-complement $K_X+D'$ of $K_X+D$ near the
fiber over $P$ there exists at most one (not necessarily
exceptional) divisor $E$ such that $a(E,D')=-1$.
\item Assume that $Z$ is a point (global case). Then $(X,D)$ is said to
be \textit{exceptional} if every $\QQ$-complement of $K_X+D$ is
klt.
\end{enumerate}
\end{definition}

\begin{definition}\label{logdel}
A pair $(X,D)$ is called {\it a
log del Pezzo surface} if the following conditions are satisfied:
\begin{enumerate}
\item every coefficient of $D$ belongs to $\M$;
\item $-(K_X+D)$ is nef;
\item $K_X+D$ is lc;
\item there exists a $\QQ$-complement of $K_X+D$.
\end{enumerate}
\end{definition}

\begin{remark} The conditions written in definition \ref{logdel} are very convenient,
although they can be replaced by other ones
(see \cite[theorem 4.1, proposition 4.6]{Sh2}). In particular, from definition
\ref{logdel} it follows that
$-(K_X+D)$ is a semi-ample divisor and a complement of $K_X+D$ exists.
\end{remark}

\begin{definition}
Define
\begin{eqnarray*}
\delta(X,D)=\#\Big\{E\mid E\ \text{is an exceptional or non-exceptional divisor}\\
\text{with discrepancy}\
a(E,D)\le -\frac67   \Big\}.
\end{eqnarray*}
\end{definition}

\begin{theorem}\cite[theorem 5.1]{Sh2}
Let $(X,D)$ be an exceptional log Del Pezzo surface.
Then $\delta(X,D) \le 2$.
\end{theorem}

\begin{proposition}\label{vidc}
Let $(X\ni P,D=bC+\sum^k_{i=1}b_iB_i)$ be a germ of two-dimensional log surface.
Assume that $b\ge 6/7$, $b_i\in \M$ and
$K_X+D$ is $\frac17$-log terminal divisor. Then
$k=1$ and one of the following three possibilities holds:
\begin{enumerate}
\item $(X\ni P,D)\an (\CC^2\ni 0,b\{x+y^2=0\}+\frac12\{x=0\})$, where $b<13/14$.
\item $(X\ni P,D)\an (\CC^2\ni 0,b\{x=0\}+b_1\{y=0\})$, where $b+b_1<13/7$.
\item $(X\ni P,D)\an (\CC^2\ni 0,b\{x=0\}+b_1\{y=0\})/\ZZ_n(q,1)$, where
$\frac{n/7-1+b_1}{1-b}<q\le n-1$ and $(q,n)=1$.
\end{enumerate}
\begin{proof} By proposition \cite[5.2]{Sh2} the pair
$(X\ni P,C+\sum^k_{i=1}b_iB_i)$ is plt, except the case (1). If $(X\ni P)$ is a
smooth point then we have the case (2).
Let $(X\ni P)$ be a singular point.
Then $(X\ni P, C)\an (\CC^2\ni 0, \{x=0\})/\ZZ_n(q,1)$ by the classification
of two-dimensional singularities \cite[theorem 2.1.3]{PrLect}.
Let us consider a cyclic cover $\psi \colon \CC^2 \to \CC^2/\ZZ_n(q,1)$ of degree $n$.
By the case (2)
$k=1$ and $\psi^{-1}(C)$ has the simple normal crossings
with $\psi^{-1}(B_1)$. Therefore we have the case (3).
\par
Take a weighted blow-up with weights  $(\frac qn,\frac 1n)$. Let $E$ be
a corresponding exceptional divisor. The condition $a(E,D)>-6/7$
is written in the case (3).
Let us prove that it is sufficient.
Let $\varphi\colon (Y,E')\to (X\ni P)$ be another blow-up
with $a(E',D)\le -6/7$. We may assume that the minimal
resolution $\widetilde X\to (X\ni P)$ factors through
$\varphi$, i.e. we have $\widetilde X\to Y\to X$.
Take the weighted blow-up of cyclic singularity $\frac 1{n_1}(q_1,1)$ at the point
$P=E'\cap C_Y$ with weights
$(\frac {q_1}{n_1},\frac 1{n_1})$, where $C_Y$ is a proper transform of $C$ on
$Y$. Then
$a(\ \cdot\ ,a(E',D)E'+bC_Y)\le -6/7$. Repeating the process we get
$a(E,D)\le -6/7$, a contradiction.
\end{proof}
\end{proposition}

Let us recall well-known statements.
\begin{proposition}\label{free2}
Let $X$ be a normal projective surface with $h^1(X,\OO_X)=0$ and $H=\sum h_iH_i$ be an
effective Cartier divisor such that $\Supp H\cap \Sing X=\emptyset$.
Assume that $p_a(H)=0$, $p_a(H_i)=0$ for all $i$ and
$H$ is nef. Then
\begin{enumerate}
\item $|H|$ is a free linear system and $\dim|H|=H^2+1$;
\item if $H^2\ge 1$ then $|H|$ gives a birational morphism
$\psi \colon X\to X'$. Besides,
$(X',\psi(H))$ is one of the following pairs: $(\PP^2,\OO_{\PP^2}(m))$, where $m=1,2$;
$(\FFF_n,E_0+mf)$, where $m\ge 1$; $(\PP(n,1,1),X_n)$, where $n\ge 2$.
\end{enumerate}
\begin{proof} From an exact sequence
$$0\longrightarrow \OO_X \longrightarrow \OO_X(H) \longrightarrow
\OO_H(H)\longrightarrow 0$$
we have
$$0\longrightarrow \CC \longrightarrow H^0(X,\OO_X(H))\longrightarrow
H^0(H,\OO_H(H))\longrightarrow 0.$$
Since $H\cdot H_i\ge 0$, where $H_i$ is an arbitrary component of
$H$ then a linear system $|H|_H$ is free \cite{CFHR}.
Therefore $|H|$ is a free linear system and
$\dim|H|=H^2+1$. Hence, it follows in the standard way that
$|H|$ separates the points if $H^2\ge 1$. Besides,
$X'$ is a minimal degree surface.
All surfaces of minimal degree are written in the statement of proposition
(see \cite[chapter 4 \S 3]{GH}).
\end{proof}
\end{proposition}

\begin{proposition}\label{free1}
Let $X$ be a normal projective surface with $h^1(X,\OO_X)=0$ and $H$ be an
irreducible curve with $p_a(H)=1$. Assume that $H^2\ge 3$ and
$\Supp H\cap \Sing X=\emptyset$. Then a linear system $|H|$ gives a birational
morphism and $\dim|H|=H^2$.
\begin{proof} The proof is the same one as in proposition \ref{free2}.
\end{proof}
\end{proposition}

\section{\bf Classification and its corollaries}
In the next theorem the classification of exceptional log Del Pezzo surfaces with
$\delta=1$ is given.
\begin{theorem}\label{mainth}
Let $(X,D_X=\sum d_iD_i)$ be an exceptional log Del Pezzo surface with
$\delta(X,D_X)=1$. If $d_i<6/7$ for all $i$ then let us take a blow-up
$\tau\colon X'\to X$ and an exceptional curve
$E$ with discrepancy $a(E,D_X)\le -6/7$ appears $(\Exc \tau=E)$.
Write $K_{X'}+D_{X'}=\tau^*(K_X+D_X)$.
If $d_i\ge 6/7$ for some $i$ then put $X'=X$.
\\
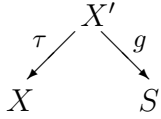
\begin{center}
\begin{picture}(80,46)(0,0)
\put(14,10){$X$}
\put(40,40){\vector(-1,-1){18}}
\put(42,42){$X'$}
\put(50,40){\vector(1,-1){18}}
\put(64,10){$S$}
\put(24,34){\footnotesize{$\tau$}}
\put(62,34){\footnotesize{$g$}}
\end{picture}
\end{center}
Then there exists a birational contraction
$g\colon X'\to S$ with the following properties: $\rho(S)=1$, $g$ doesn't contract
the curve $E$. Let $D=g(D_{X'})$. Then
$(S,D)$ is an exceptional log Del Pezzo surface with same complements
(i.e. $n$-complement of $K_X+D_X$ induces
$n$-complement of $K_S+D$ and conversely). Put $D=tC+\sum b_iB_i$, where
$6/7\le t<1$, $b_i\in \SM$, $C=g(E)$.
Then $(S,D)$ is one of the following pairs\footnote{See remark  \ref{main} about
the conditions on $C$ and $B_i$.}:
\begin{enumerate}
\item[{\bf 1)}] $S=\PP^2.$ There are the following possibilities for $D$.
\begin{enumerate}
\item[1(+4)] $D=tX_2+\frac12X_1+\frac{k-1}kX_1$, where $k=3,4$ and
$t\in[\frac67,\frac34+\frac1{2k}]$.
\item[2(+1)] $D=tX_1+\frac12X_1+\frac34X_1+\frac{k-1}kX_1$, where $k=5,6$ and
$t\in[\frac67,\frac34+\frac1k]$.
\item[3(+1)] $D=tX_1+\frac12X_1+\frac45X_1+\frac45X_1$ and
$t\in[\frac67,\frac9{10}]$.
\item[4(+1)] $D=tX_1+\frac23X_2+\frac{k-1}kX_1$, where $k=4,5$ and
$t\in[\frac67,\frac23+\frac1k]$. If $k=5$ then a straight line $X_1$
is not tangent $X_2$.
It is possible that there is $X_1+X_1$ instead of $X_2$.
\end{enumerate}

\item[{\bf 2)}] $S=\PP(1,1,2)$. There are the following possibilities for $D$.
\begin{enumerate}
\item[1({\rm ell})] $D=tX_4+\frac12X_1$ and $t\in[\frac67,\frac78]$.
\item[2(+4)] $D=tX_3+\frac23X_2$ and $t\in[\frac67,\frac89]$.
\item[3(+2)] $D=tX_2+\frac56X_2+\frac12X_1$ and $t\in[\frac67,\frac{11}{12}]$.
\item[4(+2)] $D=tX_2+\frac45X_2+\frac{k-1}kX_1$, where $k=2,3$ and
$t\in[\frac67,\frac65-\frac{k-1}{2k}]$.
\item[5(+2)] $D=tX_2+\frac34X_2+\frac{k-1}kX_1$, where $k=3,4$ and
$t\in[\frac67,\frac54-\frac{k-1}{2k}]$.
\item[6(+2)] $D=tX_2+\frac23X_2+\frac{k-1}kX_1$, where $k=4,5,6$ and
$t\in[\frac67,\frac43-\frac{k-1}{2k}]$.
\item[7(+2)] $D=tX_2+\frac12X_3+\frac{k-1}kX_1$, where $k=3,4$ and
$t\in[\frac67,\frac54-\frac{k-1}{2k}]$.
It is possible that there is $X_2+X_1$ instead of $X_3$.
\item[8(0)] $D=tX_1+\frac12X_2+\frac23X_2+\frac{k-1}kX_1$, where $k=4,5$ and
$t\in[\frac67,\frac23+\frac1k]$.
\item[9(0)] $D=tX_1+\frac12X_3+\frac45X_2$ and
$t\in[\frac67,\frac9{10}]$. It is possible that there is $X_2+X_1$ instead of $X_3$.
Note that $X_1$ and $X_3$ have the intersections at two different points.
\item[10(0)] $D=tX_1+\frac34X_2+\frac45X_2$ and $t\in[\frac67,\frac9{10}]$.
\end{enumerate}

\item[{\bf 3)}] $S=\PP(1,1,3)$. There are the following possibilities for $D$.
\begin{enumerate}
\item[1(+5)] $D=tX_4+\frac12X_3$ and $t\in[\frac67,\frac78]$.
\item[2(+3)] $D=tX_3+\frac{k-1}kX_3$, where $k=4,5$ and
$t\in[\frac67,\frac53-\frac{k-1}k]$.
\item[3(+3)] $D=tX_3+\frac12X_3+\frac{k-1}kX_1$, where $k=3,4,5,6$ and
$t\in[\frac67,\frac56+\frac1{3k}]$.
\end{enumerate}

\item[{\bf 4)}] $(+4)\ \ S=\PP(1,1,4)$ and $D=tX_4+\frac12X_5$ or
$D=tX_4+\frac12X_4+\frac12X_1$, where
$t\in[\frac67,\frac78]$.

\item[{\bf 5)}] $(+5)\ \ S=\PP(1,1,5)$, $D=tX_5+\frac12X_5$ and
$t\in[\frac67,\frac9{10}]$.

\item[{\bf 6)}] $S=\PP(1,2,3)$. There are the following possibilities for $D$.
\begin{enumerate}
\item[1({\rm ell})] $D=tX_6+\frac{k-1}kX_1$,
$t\in[\frac67,1-\frac{k-1}{6k}]$ and $2\le k\le 6$.
\item[2({\rm ell})] $D=\frac67X_7$.
\item[3(+3)] $D=tX_5+\frac12X_3$ and $t\in[\frac67,\frac9{10}]$.
\item[4(+2)] $D=tX_4+\frac{k-1}kX_3$, where $k=4,5,6$ and
$t\in[\frac67,\frac34+\frac3{4k}]$.
\item[5(+2)] $D=tX_4+\frac12X_5$ and $t\in[\frac67,\frac78]$.
It is possible that there is $X_3+X_2$ instead of $X_5$.
\item[6(+1)] $D=tX_3+\frac{k_1-1}{k_1}X_3+\frac{k_2-1}{k_2}X_2$, where
$(k_1,k_2)=(2,5)$, $(2,6)$, $(3,3)$, $(4,2)$, $(5,2)$ and
$t\in[\frac67,\frac13+\frac1{k_1}+\frac2{3k_2}]$.
\item[7(+1)] $D=tX_3+\frac12X_4+\frac23X_2$ and $t\in[\frac67,\frac89]$.
\item[8(+1)] $D=tX_3+\frac23X_5$ and $t\in[\frac67,\frac89]$.
\item[9(+1)] $D=tX_3+\frac23X_4+\frac12X_1$ and $t\in[\frac67,\frac{17}{18}]$.
\item[10(+1)] $D=tX_3+\frac{k-1}kX_4$, where $k=5,6$ and
$t\in[\frac67,\frac23+\frac4{3k}]$.
\item[11(0)] $D=tX_2+\frac12X_3+\frac23X_4$ and $t\in[\frac67,\frac{11}{12}]$.
\item[12(0)] $D=tX_2+\frac34X_3+\frac12X_4$ and $t\in[\frac67,\frac78]$.
\item[13(0)] $D=tX_2+\frac23X_3+\frac34X_3$ and $t\in[\frac67,\frac78]$.
\end{enumerate}

\item[{\bf 7)}] $S=\PP(1,2,5).$ There are the following possibilities for $D$.
\begin{enumerate}
\item[1(+3)] $D=tX_6+\frac12X_5$ and $t\in[\frac67,\frac{11}{12}]$.
\item[2(+2)] $D=tX_5+\frac23X_5$ and $t\in[\frac67,\frac{14}{15}]$.
\item[3(+2)] $D=tX_5+\frac12X_7$ and $t\in[\frac67,\frac9{10}]$.
It is possible that there is $X_5+X_2$ instead of $X_7$.
\item[4(+2)] $D=tX_5+\frac12X_6+\frac12X_1$ and $t\in[\frac67,\frac9{10}]$.
\item[5(0)] $D=tX_2+\frac12X_5+\frac34X_5$ and $t\in[\frac67,\frac78]$.
\end{enumerate}

\item[{\bf 8)}] $S=\PP(1,3,4)$. There are the following possibilities for $D$.
\begin{enumerate}
\item[1({\rm ell})] $D=tX_9$ and $t\in[\frac67,\frac89]$.
\item[2(+3)] $D=\frac67X_7+\frac12X_4$.
\item[3(+1)] $D=tX_4+\frac{k_1-1}{k_1}X_4+\frac{k_2-1}{k_2}X_3$, where
$(k_1,k_2)=(2,4)$, $(2,5)$, $(2,6)$, $(3,2)$;
$t\in[\frac67,\frac34+\frac3{4k_2}]$ if $k_1=2$, and
$t\in[\frac67,\frac{19}{21})$ if $k_1=3$.
\item[4(+1)] $D=tX_4+\frac12X_9$ and $t\in[\frac67,\frac78]$.
It is possible that there is $X_6+X_3$ instead of $X_9$.
\item[5(+1)] $D=tX_4+\frac34X_6$ and $t\in[\frac67,\frac78]$.
\item[6(0)] $D=tX_3+\frac12X_4+\frac{k-1}kX_4$, where $k=5,6$ and
$t\in[\frac67,\frac23+\frac4{3k}]$.
\item[7(0)] $D=tX_3+\frac23X_8$ and $t\in[\frac67,\frac89]$.
It is possible that there is $X_4+X_4$ instead of $X_8$.
\end{enumerate}

\item[{\bf 9)}] $S=\PP(1,3,5)$. There are the following possibilities for $D$.
\begin{enumerate}
\item[1({\rm ell})] $D=tX_{10}$ and $t\in[\frac67,\frac9{10}]$.
\item[2(+2)] $D=tX_6+\frac{k-1}kX_5$, where $k=3,4$ and
$t\in[\frac67,\frac23+\frac5{6k}]$.
\item[3(+1)] $D=tX_5+\frac12X_9$ and $t\in[\frac67,\frac9{10}]$.
It is possible that the monomial $x_1x_2x_3$ is absent in the polynomial defining
 $X_9$. Also it is possible that there is $X_6+X_3$ instead of $X_9$.
\item[4(+1)] $D=tX_5+\frac34X_6$ and $t\in[\frac67,\frac9{10}]$.
\item[5(0)] $D=tX_3+\frac23X_5+\frac12X_6$ and $t\in[\frac67,\frac89]$.
\item[6(0)] $D=tX_3+\frac12X_5+\frac34X_5$ and $t\in[\frac67,\frac{11}{12}]$.
\end{enumerate}

\item[{\bf 10)}] $S=\PP(1,2,7)$. There are the following possibilities for $D$.
\begin{enumerate}
\item[1(+3)] $D=\frac67X_7+\frac12X_8$.
\item[2(+3)] $D=tX_7+\frac12X_7$ and $t\in[\frac67,\frac{13}{14}]$.
\item[3(0)] $D=tX_2+\frac12X_7+\frac23X_7$ and $t\in[\frac67,\frac{11}{12}]$.
It is possible that $(B_1\cdot B_2)_{(0:1:0)}=3/2.$
\end{enumerate}

\item[{\bf 11)}] $S=\PP(1,4,5)$. There are the following possibilities for $D$.
\begin{enumerate}
\item[1(+1)] $D=tX_5+\frac12X_5+\frac{k-1}kX_4$, where
$k=3,4,5$;
$t\in[\frac67,\frac{19}{21})$ if $k=3$, and
$t\in[\frac67,\frac7{10}+\frac4{5k}]$ if $k\ge 4$.
\item[2(+1)] $D=tX_5+\frac23X_8$ and $t\in[\frac67,\frac{14}{15}]$.
\item[3(0)] $D=tX_4+\frac12X_5+\frac{k-1}kX_5$, where $k=4,5$ and
$t\in[\frac67,\frac58+\frac5{4k}]$.
\end{enumerate}

\item[{\bf 12)}] $S=\PP(2,3,5).$ There are the following possibilities for $D$.
\begin{enumerate}
\item[1(+1)] $D=tX_8+\frac12X_5$ and $t\in[\frac67,\frac{15}{16}]$.
\item[2(0)] $D=tX_5+\frac12X_{11}$ and $t\in[\frac67,\frac9{10}]$.
It is possible that there is $X_6+X_5$ instead of $X_{11}$.
\item[3(0)] $D=tX_5+\frac23X_8$ and $t\in[\frac67,\frac{19}{21})$.
\end{enumerate}

\item[{\bf 13)}] $S=\PP(1,3,7)$. There are the following possibilities for $D$.
\begin{enumerate}
\item[1(+2)] $D=tX_7+\frac23X_7$ and $t\in[\frac67,\frac{19}{21})$.
\item[2(+2)] $D=\frac67X_7+\frac12X_{10}$.
It is possible that there is $X_7+X_3$ instead of $X_{10}$.
\item[3(+2)] $D=tX_7+\frac12X_9$ and $t\in[\frac67,\frac{13}{14}]$.
\item[4(0)] $D=tX_3+\frac12X_7+\frac23X_7$ and $t\in[\frac67,\frac{17}{18}]$.
\end{enumerate}

\item[{\bf 14)}] $S=\PP(1,3,8).$ There are the following possibilities for $D$.
\begin{enumerate}
\item[1(+3)] $D=tX_9+\frac12X_8$ and $t\in[\frac67,\frac89]$.
\item[2(+2)] $D=tX_8+\frac12X_9$ and $t\in[\frac67,\frac{15}{16}]$.
\item[3(0)] $D=tX_3+\frac12X_8+\frac23X_8$ and $t\in[\frac67,\frac89]$.
\end{enumerate}

\item[{\bf 15)}] $S=\PP(1,4,7).$ There are the following possibilities for $D$.
\begin{enumerate}
\item[1(+2)] $D=tX_8+\frac23X_7$ and $t\in[\frac67,\frac{11}{12}]$.
\item[2(+1)] $D=\frac67X_7+\frac12X_{12}$.
It is possible that there is $X_8+X_4$ instead of $X_{12}$.
\item[3(+1)] $D=tX_7+\frac{k-1}kX_8$, where $k=3,4$ and
$t\in[\frac67,\frac47+\frac8{7k}]$.
\end{enumerate}

\item[{\bf 16)}] $(0)\ \  S=\PP(1,5,6)$, $D=tX_5+\frac12X_6+\frac34X_6$ and
$t\in[\frac67,\frac9{10}]$.

\item[{\bf 17)}] $S=\PP(2,3,7).$ There are the following possibilities for $D$.
\begin{enumerate}
\item[1(+1)] $D=tX_9+\frac12X_7$ and $t\in[\frac67,\frac{17}{18}]$.
\item[2(0)] $D=\frac67X_7+\frac12X_{12}$.
\item[3(0)] $D=\frac67X_7+\frac23X_9$.
\item[4($-1$)] $D=tX_3+\frac23X_{14}$ and $t\in[\frac67,\frac89]$.
\end{enumerate}

\item[{\bf 18)}] $S=\PP(3,4,5).$ There are the following possibilities for $D$.
\begin{enumerate}
\item[1(+1)] $D=tX_{13}$ and $t\in[\frac67,\frac{12}{13}]$.
\item[2(0)] $D=tX_9+\frac12X_8$ and $t\in[\frac67,\frac89]$.
\item[3(0)] $D=tX_8+\frac12X_9$ and $t\in[\frac67,\frac{25}{28})$.
\item[4(0)] $D=tX_8+\frac12X_{10}$ and $t\in[\frac67,\frac78]$.
\item[5($-1$)] $D=tX_5+\frac12X_{15}$ and $t\in[\frac67,\frac9{10}]$.
\end{enumerate}

\item[{\bf 19)}] $({\rm ell})\ \ S=\PP(1,5,7)$,
$D=tX_{15}$ and $t\in[\frac67,\frac{13}{15}]$.

\item[{\bf 20)}] $(+3)\ \ S=\PP(1,3,10)$,
$D=tX_{10}+\frac12X_{10}$ and $t\in[\frac67,\frac9{10}]$.

\item[{\bf 21)}] $S=\PP(1,4,9)$. There are the following possibilities for $D$.
\begin{enumerate}
\item[1(+2)] $D=tX_9+\frac12X_{12}$ and $t\in[\frac67,\frac89]$.
\item[2(0)] $D=tX_4+\frac12X_9+\frac23X_9$ and $t\in[\frac67,\frac78]$.
\end{enumerate}

\item[{\bf 22)}] $S=\PP(1,5,8)$. There are the following possibilities for $D$.
\begin{enumerate}
\item[1({\rm ell})] $D=tX_{16}$ and $t\in[\frac67,\frac78]$.
\item[2(+1)] $D=tX_8+\frac23X_{10}$ and $t\in[\frac67,\frac{11}{12}]$.
\end{enumerate}

\item[{\bf 23)}] $(0)\ \ S=\PP(3,4,7)$,
$D=tX_7+\frac12X_{15}$ and $t\in[\frac67,\frac{13}{14})$.

\item[{\bf 24)}] $(+1)\ \ S=\PP(1,5,9)$,
 $D=tX_9+\frac23X_{10}$ and $t\in[\frac67,\frac{25}{27}]$.

\item[{\bf 25)}] $S=\PP(3,5,7).$ There are the following possibilities for $D$.
\begin{enumerate}
\item[1(+1)] $D=tX_{17}$ and $t\in[\frac67,\frac{15}{17}]$.
\item[2($-1$)] $D=tX_5+\frac12 X_{21}$ and $t\in[\frac67,\frac9{10}]$.
\end{enumerate}

\item[{\bf 26)}] $(+2)\ \ S=\PP(1,4,11)$,
$D=tX_{11}+\frac12X_{12}$ and $t\in[\frac67,\frac{10}{11}]$.

\item[{\bf 27)}] $(+1)\ \ S=\PP(2,3,11)$,
$D=tX_{11}+\frac12X_{11}$ and $t\in[\frac67,\frac{21}{22}]$.

\item[{\bf 28)}] $(0)\ \ S=\PP(2,5,9)$,
$D=tX_9+\frac12X_{15}$ and $t\in[\frac67,\frac{17}{18}]$.

\item[{\bf 29)}] $(+3)\ \ S=\PP(1,4,13)$,
$D=tX_{13}+\frac12X_{13}$ and $t\in[\frac67,\frac{23}{26}]$.

\item[{\bf 30)}] $(+1)\ \ S=\PP(1,6,11)$,
$D=tX_{11}+\frac23X_{12}$ and $t\in[\frac67,\frac{10}{11}]$.

\item[{\bf 31)}] $(0)\ \ S=\PP(2,5,11)$,
$D=tX_{11}+\frac12X_{15}$ and $t\in[\frac67,\frac{21}{22}]$.

\item[{\bf 32)}] $(0)\ \ S=\PP(3,4,11)$,
$D=tX_{11}+\frac12X_{15}$ and $t\in[\frac67,\frac{21}{22}]$.

\item[{\bf 33)}] $(-1)\ \ S=\PP(3,7,8)$,
$D=\frac67X_7+\frac12 X_{24}$.

\item[{\bf 34)}] $(+2)\ \ S=\PP(1,5,13)$,
$D=tX_{13}+\frac12X_{15}$ and $t\in[\frac67,\frac{23}{26}]$.

\item[{\bf 35)}] $({\rm ell})\ \ S=\PP(1,7,11)$,
$D=tX_{22}$ and $t\in[\frac67,\frac{19}{22}]$.

\item[{\bf 36)}] $(+2)\ \ S=\PP(1,5,14)$,
$D=tX_{14}+\frac12X_{15}$ and $t\in[\frac67,\frac{25}{28}]$.

\item[{\bf 37)}] $(0)\ \ S=\PP(2,5,13)$,
$D=tX_{13}+\frac12X_{15}$ and $t\in[\frac67,\frac{25}{26}]$.

\item[{\bf 38)}] $(0)\ \ S=\PP(3,4,13)$,
$D=tX_{13}+\frac12X_{16}$ and $t\in[\frac67,\frac{12}{13}]$.

\item[{\bf 39)}] $(-1)\ \ S=\PP(4,5,11)$,
$D=tX_{11}+\frac12X_{20}$ and $t\in[\frac67,\frac{10}{11}]$.

\item[{\bf 40)}] $(+2)\ \ S=\PP(1,6,17)$,
$D=tX_{17}+\frac12X_{18}$ and $t\in[\frac67,\frac{15}{17}]$.

\item[{\bf 41)}] $(0)\ \ S=\PP(3,5,17)$,
$D=tX_{17}+\frac12X_{20}$ and $t\in[\frac67,\frac{15}{17}]$.

\item[{\bf 42)}] $(+1)\ \ S=\PP(3,4,19)$,
$D=tX_{19}+\frac12X_{19}$ and $t\in[\frac67,\frac{33}{38}]$.

\item[{\bf 43)}] $(+1)\ \ S=S(\AAA_1+\frac14(1,1)+\frac1{14}(9,1))$ is a toric
surface;
$D=tC+\frac12B_1$, where $C\sim B_1$ is a closure of one-dimensional orbit passing
through the first and second points and $t\in[\frac67,\frac{13}{14})$.

\item[{\bf 44)}] $(+1)\ \ S=S(\frac13(1,1)+\frac13(1,1)+\frac1{15}(11,1))$ is
a toric surface;
$D=tC+\frac12B_1$, where $C\sim B_1$
is a closure of one-dimensional orbit passing
through the first and second points
and $t\in[\frac67,\frac9{10}]$.

\item[{\bf 45)}] $(+1)\ \ S=S(\frac13(1,1)+\AAA_2+\frac19(4,1))$ is
a toric surface;
$D=tC+\frac12B_1$, where $B_1$
is a closure of one-dimensional orbit passing
through the first and second points
and $C\sim B_1+T$,
where $T$ is a closure of one-dimensional orbit passing
through the second and third points
$($on a minimal resolution of surface
$S$ a proper transform of $C$ passes through two $(-3)$ curves$)$
and $t\in[\frac67,\frac78]$.

\item[{\bf 46)}] $(+1)\ \ S=S(\frac14(1,1)+\AAA_1+\AAA_5)$ is
a toric surface;
$D=tC+\frac12B_1$. The structure of $C$ and $B_1$ is similar to $45(+1)$
and $t\in[\frac67,\frac9{10}]$.

\item[{\bf 47)}] $S=S(\AAA_1+\AAA_1+\AAA_3)$ is
a toric surface. There are the following possibilities for $D$.
\begin{enumerate}
\item[1(+1)] $D=tC+\frac23B_1$. The structure of $C$ and $B_1$ is similar to $45(+1)$
and $t\in[\frac67,\frac89]$.
\item[2(0)] $D=tC+\frac34B_1$, where $C$ is a closure of one-dimensional orbit passing
through the first and second points
and $B_1\sim 3T$,
where $T$ is a closure of one-dimensional orbit passing
through the first and third points
and $t\in[\frac67,\frac78]$.
\end{enumerate}

\item[{\bf 48)}] $(+1)\ \ S=S(\AAA_1+\AAA_1+\frac18(5,1))$ is
a toric surface;
$D=tC+\frac12B_1$, where $C$ is a closure of one-dimensional orbit passing
through the first and second points
and $B_1\sim C+T$,
where $T$ is a closure of one-dimensional orbit passing
through the second and third points and $t\in[\frac67,\frac78]$.

\item[{\bf 49)}] $(0)\ \ S=S(\AAA_1+\AAA_5+\frac1{16}(11,1))$ is
a toric surface;
$D=tC+\frac12B_1$, where $C$
is a closure of one-dimensional orbit passing
through the first and second points
and $B_1\sim 3T$,
where $T$  is a closure of one-dimensional orbit passing
through the first and third points and $t\in[\frac67,\frac{15}{16}]$.

\item[{\bf 50)}] $(0)\ \ S=S(\frac14(1,1)+\AAA_3+\frac1{16}(13,1))$ is
a toric surface;
$D=tC+\frac12B_1$, where $C$
is a closure of one-dimensional orbit passing
through the first and second points
and $B_1\sim 5T$,
where $T$ is a closure of one-dimensional orbit passing
through the first and third points and
$t\in[\frac67,\frac78]$.

\item[{\bf 51)}] The minimal resolution of $S$ is one of the following ones:
\\
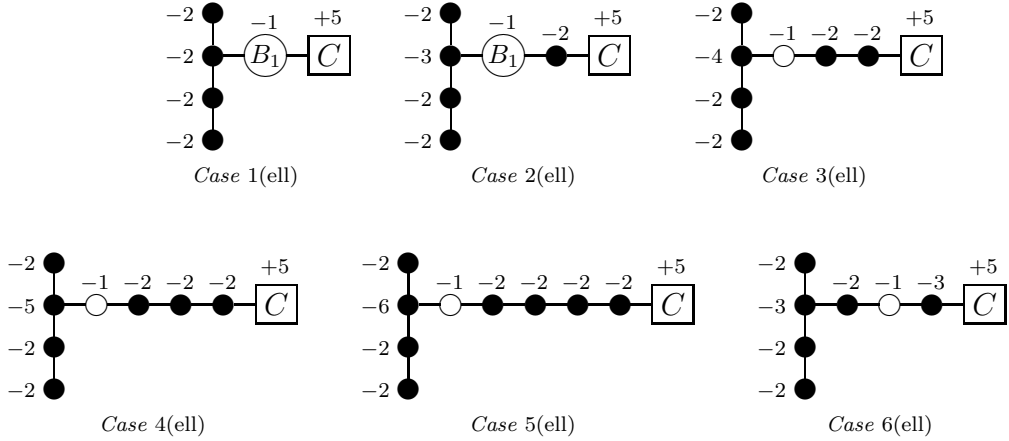
\begin{center}
\begin{picture}(290,70)(0,0)
\put(18,16){\circle*{8}}
\put(0,13){\tiny{$-2$}}
\put(18,20){\line(0,1){8}}
\put(18,32){\circle*{8}}
\put(0,29){\tiny{$-2$}}
\put(18,36){\line(0,1){8}}
\put(18,48){\circle*{8}}
\put(0,45){\tiny{$-2$}}
\put(18,52){\line(0,1){8}}
\put(18,64){\circle*{8}}
\put(0,61){\tiny{$-2$}}
\put(22,48){\line(1,0){8}}
\put(38,48){\circle{16}}
\put(32,58){\tiny{$-1$}}
\put(32,45){\footnotesize{$B_1$}}
\put(46,48){\line(1,0){8}}
\put(54,44){\fbox{$C$}}
\put(56,60){\tiny{$+5$}}
\put(10,0){\tiny{Case $1({\rm ell)}$}}
\put(108,16){\circle*{8}}
\put(90,13){\tiny{$-2$}}
\put(108,20){\line(0,1){8}}
\put(108,32){\circle*{8}}
\put(90,29){\tiny{$-2$}}
\put(108,36){\line(0,1){8}}
\put(108,48){\circle*{8}}
\put(90,45){\tiny{$-3$}}
\put(108,52){\line(0,1){8}}
\put(108,64){\circle*{8}}
\put(90,61){\tiny{$-2$}}
\put(112,48){\line(1,0){8}}
\put(128,48){\circle{16}}
\put(122,58){\tiny{$-1$}}
\put(122,45){\footnotesize{$B_1$}}
\put(136,48){\line(1,0){8}}
\put(148,48){\circle*{8}}
\put(142,54){\tiny{$-2$}}
\put(152,48){\line(1,0){8}}
\put(160,44){\fbox{$C$}}
\put(162,60){\tiny{$+5$}}
\put(115,0){\tiny{Case $2({\rm ell)}$}}
\put(218,16){\circle*{8}}
\put(200,13){\tiny{$-2$}}
\put(218,20){\line(0,1){8}}
\put(218,32){\circle*{8}}
\put(200,29){\tiny{$-2$}}
\put(218,36){\line(0,1){8}}
\put(218,48){\circle*{8}}
\put(200,45){\tiny{$-4$}}
\put(218,52){\line(0,1){8}}
\put(218,64){\circle*{8}}
\put(200,61){\tiny{$-2$}}
\put(222,48){\line(1,0){8}}
\put(234,48){\circle{8}}
\put(228,54){\tiny{$-1$}}
\put(238,48){\line(1,0){8}}
\put(250,48){\circle*{8}}
\put(244,54){\tiny{$-2$}}
\put(254,48){\line(1,0){8}}
\put(266,48){\circle*{8}}
\put(260,54){\tiny{$-2$}}
\put(270,48){\line(1,0){8}}
\put(278,44){\fbox{$C$}}
\put(280,60){\tiny{$+5$}}
\put(225,0){\tiny{Case $3({\rm ell)}$}}
\end{picture}
\end{center}

\vspace{0.5cm}

\begin{center}
\begin{picture}(320,70)(45,0)
\put(18,16){\circle*{8}}
\put(0,13){\tiny{$-2$}}
\put(18,20){\line(0,1){8}}
\put(18,32){\circle*{8}}
\put(0,29){\tiny{$-2$}}
\put(18,36){\line(0,1){8}}
\put(18,48){\circle*{8}}
\put(0,45){\tiny{$-5$}}
\put(18,52){\line(0,1){8}}
\put(18,64){\circle*{8}}
\put(0,61){\tiny{$-2$}}
\put(22,48){\line(1,0){8}}
\put(34,48){\circle{8}}
\put(28,54){\tiny{$-1$}}
\put(38,48){\line(1,0){8}}
\put(50,48){\circle*{8}}
\put(44,54){\tiny{$-2$}}
\put(54,48){\line(1,0){8}}
\put(66,48){\circle*{8}}
\put(60,54){\tiny{$-2$}}
\put(70,48){\line(1,0){8}}
\put(82,48){\circle*{8}}
\put(76,54){\tiny{$-2$}}
\put(86,48){\line(1,0){8}}
\put(94,44){\fbox{$C$}}
\put(96,60){\tiny{$+5$}}
\put(35,0){\tiny{Case $4({\rm ell)}$}}
\put(152,16){\circle*{8}}
\put(134,13){\tiny{$-2$}}
\put(152,20){\line(0,1){8}}
\put(152,32){\circle*{8}}
\put(134,29){\tiny{$-2$}}
\put(152,36){\line(0,1){8}}
\put(152,48){\circle*{8}}
\put(134,45){\tiny{$-6$}}
\put(152,52){\line(0,1){8}}
\put(152,64){\circle*{8}}
\put(134,61){\tiny{$-2$}}
\put(156,48){\line(1,0){8}}
\put(168,48){\circle{8}}
\put(162,54){\tiny{$-1$}}
\put(172,48){\line(1,0){8}}
\put(184,48){\circle*{8}}
\put(178,54){\tiny{$-2$}}
\put(188,48){\line(1,0){8}}
\put(200,48){\circle*{8}}
\put(194,54){\tiny{$-2$}}
\put(204,48){\line(1,0){8}}
\put(216,48){\circle*{8}}
\put(210,54){\tiny{$-2$}}
\put(220,48){\line(1,0){8}}
\put(232,48){\circle*{8}}
\put(226,54){\tiny{$-2$}}
\put(236,48){\line(1,0){8}}
\put(244,44){\fbox{$C$}}
\put(246,60){\tiny{$+5$}}
\put(175,0){\tiny{Case $5({\rm ell)}$}}
\put(302,16){\circle*{8}}
\put(284,13){\tiny{$-2$}}
\put(302,20){\line(0,1){8}}
\put(302,32){\circle*{8}}
\put(284,29){\tiny{$-2$}}
\put(302,36){\line(0,1){8}}
\put(302,48){\circle*{8}}
\put(284,45){\tiny{$-3$}}
\put(302,52){\line(0,1){8}}
\put(302,64){\circle*{8}}
\put(284,61){\tiny{$-2$}}
\put(306,48){\line(1,0){8}}
\put(318,48){\circle*{8}}
\put(312,54){\tiny{$-2$}}
\put(322,48){\line(1,0){8}}
\put(334,48){\circle{8}}
\put(328,54){\tiny{$-1$}}
\put(338,48){\line(1,0){8}}
\put(350,48){\circle*{8}}
\put(344,54){\tiny{$-3$}}
\put(354,48){\line(1,0){8}}
\put(362,44){\fbox{$C$}}
\put(364,60){\tiny{$+5$}}
\put(310,0){\tiny{Case $6({\rm ell)}$}}
\end{picture}
\end{center}

\vspace{0.5cm}
\begin{enumerate}
\item[1({\rm ell})] $D=tC+\frac{k-1}kB_1$, where $k=2,3$ and
$t\in[\frac67,1-\frac{k-1}{5k}]$.
\item[2({\rm ell})] $D=tC+\frac{k-1}kB_1$, where $k=1,2$ and
$t\in[\frac67,\frac{10}{11}-\frac{k-1}{11k}]$.
\item[3({\rm ell})] $D=tC$, where $t\in[\frac67,\frac{15}{17}]$.
\item[4({\rm ell})] $D=tC$, where $t\in[\frac67,\frac{20}{23}]$.
\item[5({\rm ell})] $D=tC$, where $t\in[\frac67,\frac{25}{29}]$.
\item[6({\rm ell})] $D=tC$, where $t\in[\frac67,\frac78]$.
\end{enumerate}

\item[{\bf 52)}] The minimal resolution of toric surface $S$
is one of the following ones:
\\
\begin{center}
\begin{picture}(250,95)(0,0)
\put(20,35){\circle*{8}}
\put(14,21){\tiny{$-2$}}
\put(20,39){\line(0,1){8}}
\put(20,55){\circle{16}}
\put(14,52){\footnotesize{$B_1$}}
\put(0,52){\tiny{$-1$}}
\put(20,63){\line(0,1){8}}
\put(20,75){\circle*{8}}
\put(0,72){\tiny{$-2$}}
\put(28,55){\line(1,0){8}}
\put(36,51){\fbox{$C$}}
\put(39,66){\tiny{$+4$}}
\put(52,55){\line(1,0){8}}
\put(68,55){\circle{16}}
\put(24,35){\line(1,0){16}}
\put(44,35){\circle*{8}}
\put(38,21){\tiny{$-2$}}
\put(48,35){\line(1,0){16}}
\put(68,35){\circle*{8}}
\put(62,21){\tiny{$-2$}}
\put(68,39){\line(0,1){8}}
\put(68,63){\line(0,1){8}}
\put(68,75){\circle*{8}}
\put(79,52){\tiny{$-1$}}
\put(79,72){\tiny{$-2$}}
\put(25,5){\tiny{Case $1({\rm ell)}$}}
\put(170,35){\circle*{8}}
\put(164,21){\tiny{$-3$}}
\put(170,39){\line(0,1){8}}
\put(170,55){\circle{16}}
\put(150,52){\tiny{$-1$}}
\put(170,63){\line(0,1){8}}
\put(170,75){\circle*{8}}
\put(150,72){\tiny{$-2$}}
\put(170,75){\line(1,-1){16}}
\put(186,51){\fbox{$C$}}
\put(189,66){\tiny{$+4$}}
\put(202,55){\line(1,0){8}}
\put(218,55){\circle{16}}
\put(174,35){\line(1,0){16}}
\put(194,35){\circle*{8}}
\put(188,21){\tiny{$-2$}}
\put(198,35){\line(1,0){16}}
\put(218,35){\circle*{8}}
\put(212,21){\tiny{$-2$}}
\put(218,39){\line(0,1){8}}
\put(218,63){\line(0,1){8}}
\put(218,75){\circle*{8}}
\put(229,52){\tiny{$-1$}}
\put(229,72){\tiny{$-2$}}
\put(175,5){\tiny{Case $2({\rm ell)}$}}
\put(170,79){\line(0,1){8}}
\put(170,91){\circle*{8}}
\put(150,88){\tiny{$-2$}}
\end{picture}
\end{center}

\begin{enumerate}
\item[1({\rm ell})] $D=tC+\frac12B_1$, where $t\in[\frac67,\frac78]$.
\item[2({\rm ell})] $D=\frac67C$. Note that $S=\PP(2,3,7)$.
\end{enumerate}

\item[{\bf 53)}] The minimal resolution of $S$ is one of the following ones:
\\
\begin{center}
\begin{picture}(200,70)(0,0)
\put(18,16){\circle*{8}}
\put(0,13){\tiny{$-2$}}
\put(18,20){\line(0,1){8}}
\put(18,32){\circle*{8}}
\put(0,29){\tiny{$-2$}}
\put(18,36){\line(0,1){8}}
\put(18,48){\circle*{8}}
\put(0,45){\tiny{$-2$}}
\put(18,52){\line(0,1){8}}
\put(18,64){\circle*{8}}
\put(0,61){\tiny{$-2$}}
\put(22,48){\line(1,0){8}}
\put(34,48){\circle*{8}}
\put(28,54){\scriptsize{$-2$}}
\put(38,48){\line(1,0){8}}
\put(54,48){\circle{16}}
\put(48,45){\footnotesize{$B_1$}}
\put(48,58){\tiny{$-1$}}
\put(62,48){\line(1,0){8}}
\put(70,44){\fbox{$C$}}
\put(72,60){\scriptsize{$+4$}}
\put(25,0){\tiny{Case $1({\rm ell)}$}}
\put(148,16){\circle*{8}}
\put(130,13){\tiny{$-2$}}
\put(148,20){\line(0,1){8}}
\put(148,32){\circle*{8}}
\put(130,29){\tiny{$-2$}}
\put(148,36){\line(0,1){8}}
\put(148,48){\circle*{8}}
\put(130,45){\tiny{$-2$}}
\put(148,52){\line(0,1){8}}
\put(148,64){\circle*{8}}
\put(130,61){\tiny{$-2$}}
\put(152,48){\line(1,0){8}}
\put(164,48){\circle*{8}}
\put(158,54){\scriptsize{$-3$}}
\put(168,48){\line(1,0){8}}
\put(180,48){\circle{8}}
\put(174,54){\tiny{$-1$}}
\put(184,48){\line(1,0){8}}
\put(196,48){\circle*{8}}
\put(190,54){\tiny{$-2$}}
\put(200,48){\line(1,0){8}}
\put(208,44){\fbox{$C$}}
\put(210,60){\scriptsize{$+4$}}
\put(160,0){\tiny{Case $2({\rm ell)}$}}
\end{picture}
\end{center}

\vspace{0.3cm}
\begin{enumerate}
\item[1({\rm ell})] $D=tC+\frac12B_1$, where $t\in[\frac67,\frac78]$.
\item[2({\rm ell})] $D=tC$, where $t\in[\frac67,\frac89]$.
\end{enumerate}

\item[{\bf 54)}] The minimal resolution of $S$ is the following one:
\\
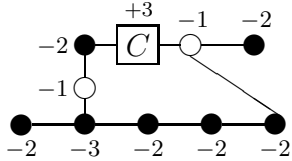
\begin{center}
\begin{picture}(120,50)(0,0)
\put(10,12){\circle*{8}}
\put(4,0){\scriptsize{$-2$}}
\put(14,12){\line(1,0){16}}
\put(34,12){\circle*{8}}
\put(28,0){\scriptsize{$-3$}}
\put(38,12){\line(1,0){16}}
\put(58,12){\circle*{8}}
\put(52,0){\scriptsize{$-2$}}
\put(62,12){\line(1,0){16}}
\put(82,12){\circle*{8}}
\put(76,0){\scriptsize{$-2$}}
\put(86,12){\line(1,0){16}}
\put(106,12){\circle*{8}}
\put(100,0){\scriptsize{$-2$}}
\put(34,14){\line(0,1){8}}
\put(34,26){\circle{8}}
\put(16,23){\scriptsize{$-1$}}
\put(34,30){\line(0,1){8}}
\put(34,42){\circle*{8}}
\put(16,39){\scriptsize{$-2$}}
\put(38,42){\line(1,0){8}}
\put(46,38){\fbox{$C$}}
\put(49,53){\tiny{$+3$}}
\put(62,42){\line(1,0){8}}
\put(74,42){\circle{8}}
\put(69,49){\scriptsize{$-1$}}
\put(74,38){\line(3,-2){34}}
\put(78,42){\line(1,0){16}}
\put(98,42){\circle*{8}}
\put(93,49){\scriptsize{$-2$}}
\end{picture}
\end{center}

$({\rm ell})\ D=\frac67C$.

\item[{\bf 55)}] The minimal resolution of $S$ is the following one:
\\
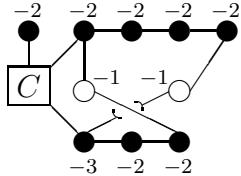
\begin{center}
\begin{picture}(85,80)(0,0)
\put(10,40){\fbox{$C$}}
\put(18,51){\line(0,1){10}}
\put(18,65){\circle*{8}}
\put(12,70){\tiny{$-2$}}
\put(26,51){\line(1,1){11}}
\put(39,65){\circle*{8}}
\put(33,70){\tiny{$-2$}}
\put(43,65){\line(1,0){10}}
\put(57,65){\circle*{8}}
\put(51,70){\tiny{$-2$}}
\put(61,65){\line(1,0){10}}
\put(75,65){\circle*{8}}
\put(69,70){\tiny{$-2$}}
\put(79,65){\line(1,0){10}}
\put(93,65){\circle*{8}}
\put(87,70){\tiny{$-2$}}
\put(26,37){\line(1,-1){11}}
\put(39,23){\circle*{8}}
\put(33,11){\tiny{$-3$}}
\put(43,23){\line(1,0){10}}
\put(57,23){\circle*{8}}
\put(51,11){\tiny{$-2$}}
\put(61,23){\line(1,0){10}}
\put(75,23){\circle*{8}}
\put(69,11){\tiny{$-2$}}
\put(39,42){\circle{8}}
\put(39,61){\line(0,-1){15}}
\put(42,45){\tiny{$-1$}}
\put(43,42){\line(2,-1){31}}
\put(75,42){\circle{8}}
\put(60,45){\tiny{$-1$}}
\put(71,42){\line(-2,-1){10}}
\put(58,35){\oval(6,6)[tr]}
\put(51,32){\line(-2,-1){12}}
\put(53,35){\oval(6,6)[bl]}
\put(79,42){\line(2,3){14}}
\end{picture}
\end{center}

$(0)$\ $D=tC$, where $t\in[\frac67,\frac{10}{11}]$.

\item[{\bf 56)}] The minimal resolution of $S$ is the following one:
\\
\begin{center}
\begin{picture}(100,105)(0,0)
\put(10,50){\fbox{$C$}}
\put(26,54){\line(1,0){10}}
\put(40,54){\circle*{8}}
\put(34,59){\tiny{$-2$}}
\put(44,54){\line(1,0){10}}
\put(58,54){\circle*{8}}
\put(52,59){\tiny{$-2$}}
\put(62,54){\line(1,0){10}}
\put(76,54){\circle*{8}}
\put(70,59){\tiny{$-2$}}
\put(40,12){\circle*{8}}
\put(44,12){\line(1,0){10}}
\put(34,0){\tiny{$-2$}}
\put(58,12){\circle*{8}}
\put(52,0){\tiny{$-3$}}
\put(62,12){\line(1,0){10}}
\put(76,12){\circle*{8}}
\put(70,0){\tiny{$-2$}}
\put(40,13){\line(-1,2){17}}
\put(40,94){\circle*{8}}
\put(34,99){\tiny{$-2$}}
\put(44,94){\line(1,0){28}}
\put(76,94){\circle*{8}}
\put(70,99){\tiny{$-2$}}
\put(40,95){\line(-1,-2){17}}
\put(96,54){\circle{8}}
\put(96,59){\tiny{$-1$}}
\put(76,94){\line(1,-2){18}}
\put(58,12){\line(1,1){38}}
\put(40,12){\line(0,1){20}}
\put(40,36){\circle{8}}
\put(30,41){\tiny{$-1$}}
\put(76,54){\line(-2,-1){33}}
\put(40,54){\line(1,-1){10}}
\put(53,42){\oval(6,6)[tl]}
\put(65,37){\circle{8}}
\put(67,41){\tiny{$-1$}}
\put(53,43){\oval(6,6)[br]}
\put(56,41){\line(1,0){10}}
\put(76,12){\line(0,1){11}}
\put(74,26){\oval(6,6)[br]}
\put(71,25){\oval(6,6)[tl]}
\put(69,37){\line(0,-1){10}}
\end{picture}
\end{center}

$(0)$\ $D=\frac67C$.

\end{enumerate}
\end{theorem}

\begin{remark}\label{main}
In the elliptic curve case $(p_a(C)=1)$ we always suppose that the singular point of
$C$ (if it exists) is an ordinary double point and every component
$B_i$ doesn't pass through it.
\par
If $b_i=1/2$ then the intersection multiplicity of
$C$ and $B_i$ is not more then 2 in the smooth point of surface.
If $b_i\ge 2/3$ then it is equal to 1 in the smooth point of S.
Consider a singular point of S. Then
$C$ and $B_i$ are the different components of toric boundary, i.e.
$K_S+B_i+C\sim 0$ is a lc divisor in the neighborhood of singularity.
The details are given in proposition \ref{vidc}.
\par
In many cases it is enough to require the irreducibility and reducibility of $X_d$.
The reader can easily find the required conditions in every case.
The variants of
$D$ which lead to different
$\widehat D$ (see the definition of $\widehat D$ in \S 3) are shown in the theorem.
For instance, see case
$9-3(+1)$. Also the minimal complementary index of $(S,D)$ can be easily computed.
\par
If we write ell in the brackets then
$p_a(C)=1$. If we write $q$ in the brackets then $p_a(C)=0$ and the self-intersection
index of proper transform of
$C$ on a minimal resolution of
$S$ is equal to $q$.
\par
In theorem
\ref{mainth} the toric surfaces not being the weighted projective spaces are
written out.
The reader will have no difficulty in finding their minimal resolutions.
Also, since we know the numbers
$q$ and $r$, where
$r$ is the number of singularities of $S$ lying on $C$
then it is easily to find this surface in the text.
\end{remark}

\begin{corollary}
In the notations of theorem \ref{mainth} the surface $S$ is toric, except
the cases $51$, $53$, $54$, $55$, $56$.
\end{corollary}

By theorem \ref{mainth} and \cite[\S 5]{Sh2} we have the next corollary.
\begin{corollary}
Let $(X,D_X)$ be an exceptional log Del Pezzo surface with
$\delta(X,D_X)\ge 1$. Then the minimal complementary index
${\rm Compl}(X,D_X)\le 66$, where
$$
{\rm Compl}(X,D_X)=\min\{n\mid K_X+D_X\ \text{is $n$-complementary} \}.
$$
\end{corollary}

\begin{remark}
Let us remark that the case ${\rm Compl}(X,D_X)=66$ is realized. Then
$(X,D_X)=(\PP^2,\frac12L_1+\frac23L_2+\frac{10}{11}L_3+\frac{12}{13}L_4)$, where
$L_i$ are the straight lines in the general position.
Also this surface can be realized as an exceptional divisor of plt blow-up of
three-dimensional exceptional canonical hypersurface singularity
$x_1^2+x_2^3+x_3^{11}+x_4^{13}=0\subset (\CC^4,0)$. Hypothetically,
${\rm Compl}(X,D_X)\le 66$ for any log Del Pezzo surface with standard coefficients.
In particular, the same conjecture ${\rm Compl}(X,D_X)\le 66$ for the
three-dimensional contractions of local type follows from the previous
two-dimensional conjecture.
\end{remark}

\begin{definition} A log canonical threshold of $(X,D)$ at the point $P$ is
denoted by $c_P(X,D)$. For every $n\in\ZZ_{>0}$ define the set
$\TT_n\subset [0,1]$
\[
\mathcal{T}_n:= \left\{c(X,F)\ \left|\ \begin{array}{c}
\text{$\dim X=n$, $(X\ni P)$ has log canonical singularities}\\
\text{and $D\ne 0$
is an effective Weil $\QQ$-Cartier divisor}
\end{array}\right.\right\}.
\]
\end{definition}

\begin{corollary} \cite[\S 5]{Kollar}
$\mathcal T_3\cap(\frac{41}{42},1)=\emptyset$.
\begin{proof} Let $6/7\le c=c_P(X,D)<1$. Let
$\psi\colon Y\to X$ be an inductive blow-up of $(X,cD)$ \cite[theorem 1.5]{Kud2}.
Then $K_Y+E+cD_Y=\psi^*(K_X+cD)$. A pair $(E,\Diff_E(cD_Y))$ is a log Del Pezzo
surface by corollary 3.10 \cite{Sh1}.
If it is non-exceptional then there exists 1-,2-,3-, 4- or 6-complement by theorem
2.3 \cite{Sh2}, a contradiction with
$c<1$. Thus, $(X,D)$ corresponds to the unique exceptional log Del Pezzo surface
with
$\delta\ge 1$. The uniqueness follows from the uniqueness of inductive blow-up
of exceptional pair.
Also $\psi(E)=P$
(cf. corollary 1.7 and proposition 1.8 \cite{Kud2}).
Now our corollary is proved by exhaustion of all cases in
theorem \ref{mainth} and \cite[\S 5]{Sh2}.
\end{proof}
\end{corollary}

\begin{remark} Using theorem \ref{mainth} and \cite[\S 5]{Sh2} the reader will
easily describe the finite set
$\mathcal T_3\cap(\frac67,1)$.
\end{remark}

Let us consider some examples demonstrating the inductive connection of
log Del Pezzo surfaces and three-dimensional singularities
\cite{Kud1}.
The calculation details can be found in \cite{Kud1}.

\begin{example}\label{ex1}
Consider the exceptional canonical singularity\\
$(X\ni 0)=(t^2+z^3+x^7y^4+azy^8+by^{12}=0\subset (\CC^4,0))$, where $|a|+|b|\ne 0$.
It is 7-complementary \cite{Kud1}.
A weighted blow-up of $\CC^4$ with weights $(42,28,8,7)$
induces a plt blow-up
$(Y,E)\to (X\ni 0)$. Then
\begin{gather*}
\big(E,\Diff_E(0)\big)=\big(t^2+z^3+xy^2+azy^4+by^6\subset \PP(3,2,4,1),\\(6/7)\{x=0\}+
(1/2)\{y=0\}\big).
\end{gather*}

The singularities of $E$ are
$\big(x_1^2+x_2^3+x_3^2=0\subset (\CC^3,0)\big)/\ZZ_4(3,2,1)=\DDD_5$
(see \cite[point 4.10]{YPG}).
Since $K_E^2=4$ then $K_{\widetilde E}^2=4$, where $\widetilde E$ is a minimal
resolution of $E$. By Noether's formula $\rho(\widetilde E)=6$ and
$\rho(E)=1$. The curve $\{x=0\}$ is elliptic. We get $E=S(\DDD_5)$.
It is the case $53-1({\rm ell})$.
\end{example}

\begin{example}
Consider the exceptional canonical singularity\\
$(X\ni 0)=(f=t^2+z^3x+x^7y-z^2y^4=0\subset (\CC^4,0))$.
It is 9-complementary \cite{Kud1}.
By the same argument as in the previous example we have
$(E,\Diff_E(0))=(f\subset \PP(43,25,11,9),0)$. The singularities of $E$ are
$\frac1{25}(2,1)$, $\frac1{11}(1,8)$ and
$\big(x_1^2+x_2^2+x_3^7=0\subset (\CC^3,0)\big)/\ZZ_9(7,7,2)=\frac1{63}(55,1)$.
Calculating $K_E^2=\frac{8}{25\cdot 11\cdot 9}$ and
$K_{\widetilde E}^2=-11$ we get $\rho(E)=2$.
Since $\delta(E,0)=1$ then let us consider a blow-up $E'\to E$.
The unique exceptional curve $C'$ has the discrepancy
$a(C',0)=-\frac{22}{25}\le -\frac67$.
We have $\{x=0\}=\{x=t-zy^2=0\}\cup\{x=t+zy^2=0\}=C_1\cup C_2$.
The minimal resolution $\widetilde E\to E$ is shown in the next figure.
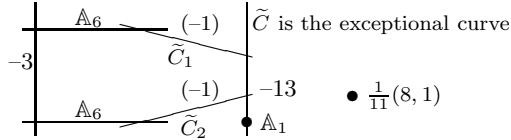
\begin{center}
\begin{picture}(140,60)(0,0)
\put(10,0){\line(0,1){50}}
\put(0,25){\scriptsize{--3}}
\put(5,5){\line(1,0){55}}
\put(25,7){\scriptsize{$\AAA_6$}}
\put(5,40){\line(1,0){55}}
\put(25,42){\scriptsize{$\AAA_6$}}
\put(90,0){\line(0,1){50}}
\put(90,5){\circle*{4}}
\put(95,2){\scriptsize{$\AAA_1$}}
\put(95,15){\scriptsize{--13}}
\put(42,42){\line(4,-1){50}}
\put(42,3){\line(4,1){50}}
\put(65,0){\tiny{$\widetilde C_2$}}
\put(65,15){\tiny{(--1)}}
\put(65,40){\tiny{(--1)}}
\put(60,27){\tiny{$\widetilde C_1$}}
\put(92,40){\tiny{$\widetilde C$}}
\put(99,40){\tiny{\ is the exceptional curve}}
\put(130,15){\circle*{4}}
\put(135,13){\tiny{$\frac1{11}(8,1)$}}
\end{picture}
\end{center}
\end{example}

Let us contract the proper transforms of $C_1$ and $C_2$ on $E'$.
We get $\psi\colon E'\to S$.
Then $p_a(\psi(C'))=1$. We obtain the surface from the case $51-2({\rm ell})$.

\begin{example}\label{ex3}
Consider the exceptional canonical singularities\\
$(X\ni 0)=(t^2+z^3+x^9y+x^4y^n=0\subset (\CC^4,0))$, where $n=7,9$.
They are 18-, 30-complementary for $n=7,9$ respectively \cite{Kud1}.
Then
\begin{gather*}
\big(E,\Diff_E(0)\big)=\big(t+z+x^9y+x^4y^n\subset \PP(9n-4,9n-4,n-1,5),\\(1/2)\{t=0\}+
(2/3)\{z=0\}\big)=\big(\PP(9n-4,n-1,5),\Diff_E(0)\big).
\end{gather*}

The singularities of $E$ are
$\frac1{n-1}(1,1)$, $\AAA_4$, $\frac1{9n-4}(2n-1,1)$ and $\delta(E,\Diff_E(0))=2$.
Let us extract the two curves with discrepancies
$\frac5{6n-6}-1\le -6/7$ and
$-\frac{14}{15}\le -6/7$ and contract the proper transforms of
$\{x=0\}$ and
$\{y=0\}$. We get $\PP^2$. It is the type $A^1_2$ \cite[\S 5]{Sh2}.
\end{example}

\section{\bf Beginning of main theorem proof}
The existence of $g\colon X'\to S$ with required properties was proved in \S 4, \S 5
of the paper \cite{Sh2}.
\par
Now we introduce the basic notions and notations used later on.
\par
Always we assume that
$(S,D=tC+\sum b_iB_i)$ is an exceptional log Del Pezzo surface with
$\delta(S,D)=1$, $6/7 \le t<1$,
$b_i\in \SM$ and $\rho(S)=1$.
It is clear that  $S$ is a rational surface.
\par
Let us define a rational number $b$ from the following equality
$K_S+D'\stackrel{\rm def}{=}K_S+bC+\sum b_iB_i\equiv 0$.
It  can happen that $\delta(S,D')=2$.

\begin{definition}
Let $g'\colon S^{\rm min}\to S$ be a minimal resolution of singularities
$P_{i_1}$,\ldots, $P_{i_m}$ lying on $C$.
Let us contract all curves from $\Exc g'$, which don't intersect a proper
transform of
$C$ on $S^{\rm min}$. We get $S^{\rm min}\to \widetilde S
\stackrel{f}{\to} S$. The surface $\widetilde S$ is called
{\it a partial resolution of} $S$
{\it along $C$ taking at the points $P_{i_1}$,\ldots, $P_{i_m}$}.
If $\Sing S\cap C=\{P_{i_1}$,\ldots, $P_{i_m}\}$ then the surface
$\widetilde S$ is called
{\it a partial resolution of} $S$
{\it along $C$}.
\par
The proper transforms of
$C$ and $B_i$ on $\widetilde S$ are denoted by
$\widetilde C$ and $\widetilde B_i$ respectively.
\end{definition}

The exceptional curves from $\Exc f$ are denoted by
$\widetilde E_1$,\ldots,$\widetilde E_r$.
It is obvious that $r$ is a number of singularities of $S$, which lie on
$C$. By proposition \ref{vidc} the singularities lying on $C$ are
$\CC^2_{x,y}/\ZZ_{n_i}(q_i,1)$, where $i=1,\ldots,r$.
The curve $C$ is given by the equation $x=0$ at the points $P_i$ .
\\
\begin{center}
\begin{picture}(80,50)(0,0)
\put(14,15){$S$}
\put(40,45){\vector(-1,-1){16}}
\put(42,47){$\widetilde S$}
\put(50,45){\vector(1,-1){16}}
\put(64,15){$\widehat S$}
\put(24,39){\footnotesize{$f$}}
\put(62,39){\footnotesize{$h$}}
\put(27,0){Fig. 1}
\end{picture}
\end{center}

Thus
$$
K_{\widetilde S}+b\widetilde C+\sum b_i\widetilde B_i+\sum_{i=1}^ra_i
\widetilde E_i=f^*(K_S+bC+\sum b_iB_i).
$$
\par
Let $\widetilde S$ be a partial resolution of $S$ along $C$.
Then we will construct a birational morphism $h\colon \widetilde S\to \widehat S$
in the case
$p_a(C)=1$ and in the case $p_a(C)=0$, $\widetilde C^2\ge 1$, where
$\widehat S$ will be a well-known surface.
The morphism $h$ will be given by a linear system $|\widetilde C|$.
\par
In the case $p_a(C)=0$ and $\widetilde C^2=0$ the
birational morphism $f$ will be the composition
of partial resolutions.
The birational morphism $h$ will be given by a linear system
$|\widetilde E_1+m\widetilde C|$, where $m\gg 0$.
\par
In the case $p_a(C)=0$ and $\widetilde C^2=-1$ the surface $\widetilde S$
will be constructed by the following way:
$S\stackrel{\varphi}{\longleftarrow} S^{\circ}\stackrel{\psi}{\longrightarrow}
\widetilde S$, where $\varphi$ is the composition
of partial resolutions and
$\psi$ is the contraction of proper transform of $C$. Also a birational morphism
$h$ will be given by a concrete linear system.
\par
In the first and second possibilities the birational morphism
$h$ doesn't contract $\widetilde C$ and
$\widetilde E_i$ for all $i$. We have
\begin{gather*}
K_{\widetilde S}+b\widetilde C+\sum b_i\widetilde B_i+
\sum_{i=1}^ra_i\widetilde E_i=h^*(K_{\widehat S}+b\widehat C+
\sum b_i\widehat B_i+\sum_{i=1}^r a_i\widehat E_i)=\\=h^*(K_{\widehat S}+\widehat D).
\end{gather*}
Note that $h$ can contract some $\widetilde B_i$.
\par
It is clear that $\big(\widehat S, t\widehat C+\sum b_i\widehat B_i+
\sum_{i=1}^r a(\widetilde E_i,D)\widehat E_i\big)$ is an exceptional log surface.
It is easily to prove the next lemma by proposition \ref{vidc}.

\begin{lemma}\label{vida}
We have $a(\widetilde E_i,bC)=1-\frac{(1-b)q_i+1}{n_i}$. This implies that
$a_i\ge 3/7$ and if $a_i=3/7$ then $n_i=2$, $q_i=1$, $b=6/7$,
$P_i\notin \Supp  B_k$ for all $k$.
\end{lemma}

\par
The idea of classification is the following one
(cf. \cite{Kud4}, \cite{Kud5}): since we know the structure of
$\widehat S$ we sort out all possibilities for
$\widehat D$. Then we describe all possible birational morphisms
$h$ and $f$ for every $\widehat D$.
By the construction it is clear that any exceptional divisor
$E$ of $h$ has a discrepancy
$a(E,\widehat D)$=0 or
$(\frac1n-1)$ and $f$ must contract all proper transforms of irreducible divisors from
$\widehat D$ with non-standard coefficients, except a proper transform of $\widehat C$.

\begin{definition} The birational morphism $h$ is called {\it an extraction}.
Every curve from $\Exc h$ is called {\it an extracted curve}.
\end{definition}

\begin{proposition}\cite[proposition 5.4]{Sh2} $p_a(C)\le 1$.
\end{proposition}

In section 4 the elliptic curve case (i.e. $p_a(C)=1$) is considered.
In sections 5, 6, 7 the rational curve case (i.e. $p_a(C)=0$ and
$\widetilde C^2\ge 1$, $\widetilde C^2=0$, $\widetilde C^2= -1$ respectively)
is considered.

\section{\bf Elliptic curve case: ${\bf p_a(C)=1}$}

Let $p_a(C)=1$.
According to proposition \ref{vidc} a singularity of curve $C$ can be only
ordinary double point.

\begin{proposition}\label{CC}
$\deg\Diff_C(0)=\sum_{i=1}^r (1-1/n_i)$ and
$$C^2=\frac{\deg \Diff_C(0)+\sum b_i(B_i\cdot C)}{1-b},
\ \ \ \widetilde C^2=C^2-\sum_{i=1}^r\frac{q_i}{n_i}.$$
\begin{proof} By the adjunction theorem $\deg \Diff_C(0)=(K_S+C)\cdot C=
(K_S+bC+\sum b_iB_i)\cdot C+
((1-b)C-\sum b_iB_i)\cdot C=(1-b)C^2-\sum b_i(B_i\cdot C)$.
\end{proof}
\end{proposition}

\begin{corollary}\label{corell1}
We have $\widetilde C^2\ge 3$.
This implies that a linear system $|\widetilde C|$ gives a birational morphism
by proposition \ref{free1}.
\begin{proof} If $\Diff_C(0)\ne \emptyset$ then
$\widetilde C^2\ge C^2-\sum_{i=1}^r\frac{n_i-1}{n_i}\ge
6\cdot\sum_{i=1}^r\frac{n_i-1}{n_i}\ge 3.$
If $\Diff_C(0)=\emptyset$ then
$\widetilde C^2\ge 7\cdot\sum b_i(B_i\cdot C)\ge 7/2$ and hence
$\widetilde C^2\ge 4$.
\end{proof}
\end{corollary}

In figure 1 the birational morphism
$h$ is given by a linear system
$|\widetilde C|$.

\begin{proposition}\label{model}
We have $r\le 2$. Moreover,
$r=0$ if and only if $\widetilde C$ is very ample divisor.
\begin{proof} Let $r\ge 3$. Then as in corollary \ref{corell1} $\widetilde C^2\ge
6\cdot\sum_{i=1}^r\frac{n_i-1}{n_i}\ge 9$. Let $\overline{S}$ be a minimal model of
$S^{\rm min}$ and let $\overline{C}$, $\overline{E}_i$ be the images of
$C$, $E_i$ on $S^{\rm min}$ respectively. Since $p_a(\overline{C})=1$
\cite[proposition 5.4]{Sh2} then $\overline{S}=\PP^2,\FFF_n$, where $n=0,2$ and
$\overline{C}\sim\OO_{\PP^2}(3),\ 2E_0+2f,\ 2E_0$ respectively.
Since $\Supp \widetilde C \cap \Sing \widetilde S=\emptyset$ then
$\overline{C}^2\ge 9$. Therefore $\overline{S}=\PP^2$, $\overline{C}^2=9$
and hence the proper transforms of curves $\widetilde E_i$
don't contract in the transfer process to the minimal model.
By lemma \ref{vida} we get a contradiction with
$-(K_{\PP^2}+b\overline{C}+\sum_{i=1}^3a_i\overline{E}_i)$ to be nef.
\par
If $r=0$ then $S=\widetilde S=\widehat S$. Now, let $\widetilde C$ be very ample
divisor and $r\ge 1$. Then as above $p_a(\overline{C})=1$. Therefore
every (--1) curve intersects $\widetilde C$ only at the single point
in the transfer process to the minimal model
$\overline{S}$. Hence
$K_{S^{\rm min}}+C_{S^{\rm min}}=\psi^*(K_{\overline{S}}+\overline{C})\equiv 0$,
where $\psi\colon S^{\rm min}\to \overline{S}$. We get a contradiction
$0\ge 2+(E_i^{\rm min})^2= -K_{S^{\rm min}}\cdot E_i^{\rm min}=
C_{S^{\rm min}}\cdot E_i^{\rm min}=1$,
where $E_i^{\rm min}$ is a proper transform of $E_i$.
\end{proof}
\end{proposition}

\begin{corollary} Just one of the following two possibilities holds:
\begin{enumerate}
\item $r\le 2$ and $\rho(\widehat S)=1;$
\item $r=2$ and $\rho(\widehat S)=2.$
\end{enumerate}
\begin{proof} Since $r\le 2$ then $\rho(\widehat S)\le 3$. If
$\rho(\widehat S)=3$ then $r=2$ and hence $\widetilde C$ is very ample divisor,
a contradiction. Similarly, if $\rho(\widehat S)=2$ then
$r\ne 0,1$. Thus $r=2$.
\end{proof}
\end{corollary}

\subsection{} {\it Consider the first case $\rho(\widehat S)=1$ and $r\le 2$.}
\par
We have $(K_{\widehat S}+\widehat C)\cdot \widehat C=0$. Therefore
$-K_{\widehat S}\sim \widehat C$ is Cartier divisor, i.e. $\widehat S$ is log
Del Pezzo surface with Du Val singularities.
Since $K_{\widehat S}^2\ge 3$ then
$\widehat S=\PP^2$, $\PP(1,1,2)$, $\PP(1,2,3)$, $S(\AAA_4)$, $S(2\AAA_1+\AAA_3)$,
$S(\DDD_5)$, $S(3\AAA_2)$, $S(\AAA_1+\AAA_5)$ or $S(\EEE_6)$ \cite{Fur}.
The curves $\widehat E_i$ are (--1) curves on a minimal resolution of
$\widehat S$ because
$\widehat C\cdot \widehat E_i=1$. Therefore $r$ is not more then the number of (--1)
curves on the minimal resolution of
$\widehat S$.

\begin{lemma}\label{CC1}
Assume that $\widetilde C^2\le 6$, i.e. $\widehat S\ne \PP^2,\ \PP(1,1,2)$
and let $\widehat B_i\ne 0$.
Then $\widehat B_i\cdot \widehat C=1$ and $\widehat B_i$ is $(-1)$
curve on a minimal resolution of
$\widehat S$.
\begin{proof}
Indeed, if $\widehat B_i\cdot \widehat C\ge 2$ then
$\widetilde C^2\ge 7-\frac{n_1-1}{n_1}$ by proposition \ref{CC}. Hence
$\widetilde C^2\ge 7$, a contradiction.
\end{proof}
\end{lemma}

Let us consider case by case all variants of $\widehat S$.
\par {\bf A).} $\widehat S=\PP^2$. Then $r=0$. There are no cases.

\par {\bf B).} $\widehat S=\PP(1,1,2)$. Then $r=0$.
We have the case 2--1({\rm ell}).

\par {\bf C).} $\widehat S=\PP(1,2,3)$. Then $r=0,1$.
If $r=0$ then we have the case 6--1({\rm ell}). Let $r=1$. By lemma \ref{vida}
we have $\widehat D=b\widehat C+6(1-b)\widehat E_1$, where
$\widehat C_1=X_6$ and $\widehat E_1=X_1$. If $b=6/7$ then $6(1-b)=6/7$.
Therefore we can assume that
$b>6/7$. Our problem is reduced to describe the following procedures.
At first take a blow-up
$h\colon \widetilde S\to \widehat S$ with the single exceptional divisor
$E$ and a discrepancy $a(E,6(1-b)\widehat E_1)$ being equal to
0 or $1/n-1$. We also require that a curve
$\widetilde E_1$ has a self-intersection index
$\le -2$ on a minimal resolution of $\widetilde S$.
After it we contract
$\widetilde E_1$. The surface obtained is a required one.
In our variant the extraction of necessary curve can happen only at
two singular points $\AAA_1$ and $\AAA_2$.
\par
Consider the first opportunity $(\CC^2,6(1-b)\{x=0\})/\ZZ_2(1,1)$.
Under the condition $b>6/7$ the extraction of required curve is shown
in the next diagram.
\\
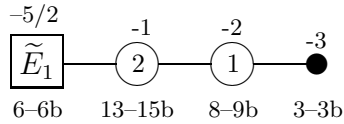
\begin{center}
\begin{picture}(150,40)(0,0)
\put(10,20){\fbox{$\widetilde E_1$}}
\put(11,4){\scriptsize{6--6b}}
\put(10,40){\scriptsize{--5/2}}
\put(30,24){\line(1,0){20}}
\put(58,24){\circle{16}}
\put(44,4){\scriptsize{13--15b}}
\put(56,21){\footnotesize{2}}
\put(56,35){\scriptsize{-1}}
\put(66,24){\line(1,0){20}}
\put(94,24){\circle{16}}
\put(85,4){\scriptsize{8--9b}}
\put(92,21){\footnotesize{1}}
\put(90,35){\scriptsize{-2}}
\put(102,24){\line(1,0){20}}
\put(126,24){\circle*{8}}
\put(117,4){\scriptsize{3--3b}}
\put(122,30){\scriptsize{-3}}
\end{picture}
\end{center}

\begin{definition} Let us describe the diagram of such type.
The numbers over circles are equal to the self-intersection indexes of exceptional
curves.
The number over square is equal to the difference of self-intersection
indexes of corresponding curve.
The numbers below circles and squares are equal to the corresponding discrepancies
taken with an opposite sign for a convenience.
The required extracted curves are enumerated.
\par
Let $k$ curves be enumerated. If we have to extract $k'$ curves then we have to
contract the remained curves
(of course we have to
contract the remained enumerated $k-k'$ curves).
Thus we have the different $C^{k'}_k$ possibilities.
\par
Unless otherwise stated we assume that the curve enumerated is extracted with
a discrepancy 0.
\end{definition}

If the curve 1 is extracted (and $b=8/9$) then we have the case 8--1({\rm ell}).
If the curve 2 is extracted (and $b=13/15$) then we have the case 19(ell).
The obtained surfaces $S$ are toric. In order to not check every time
that the extractions and contractions are the toric ones we can use
a toric criterion.

\begin{theorem}\cite[theorem 6.4]{Sh2}
A normal projective surface
$X$ is toric if and only if
there exists a boundary $D_X=\sum_{i=1}^mD_i$
such that $K_X+D_X\equiv 0$ is a lc divisor and $\rho(X)=m-2$.
\end{theorem}

In this and next cases we take a standard toric boundary $T=\sum_{i=1}^3\{x_i=0\}$
on $\widehat S$. After taking blow-up
$h$ the exceptional divisor appears with
discrepancy
$a(\ \cdot\ ,T)=-1$. After taking contraction $f$ we obtain a desired toric boundary
consisting of three required divisors.
Thus $S$ is a toric surface.
\par
Consider the second opportunity $(\CC^2,6(1-b)\{x=0\})/\ZZ_3(2,1)$. Then
\\
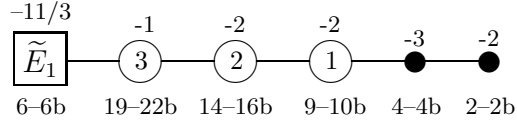
\begin{center}
\begin{picture}(196,40)(0,0)
\put(10,20){\fbox{$\widetilde E_1$}}
\put(11,4){\scriptsize{6--6b}}
\put(9,40){\scriptsize{--11/3}}
\put(30,24){\line(1,0){20}}
\put(58,24){\circle{16}}
\put(44,4){\scriptsize{19--22b}}
\put(56,21){\footnotesize{3}}
\put(56,35){\scriptsize{-1}}
\put(66,24){\line(1,0){20}}
\put(94,24){\circle{16}}
\put(80,4){\scriptsize{14--16b}}
\put(92,21){\footnotesize{2}}
\put(90,35){\scriptsize{-2}}
\put(102,24){\line(1,0){20}}
\put(130,24){\circle{16}}
\put(120,4){\scriptsize{9--10b}}
\put(128,21){\footnotesize{1}}
\put(126,35){\scriptsize{-2}}
\put(138,24){\line(1,0){20}}
\put(162,24){\circle*{8}}
\put(153,4){\scriptsize{4--4b}}
\put(158,30){\scriptsize{-3}}
\put(166,24){\line(1,0){20}}
\put(190,24){\circle*{8}}
\put(181,4){\scriptsize{2--2b}}
\put(186,30){\scriptsize{-2}}
\end{picture}
\end{center}

If the curve 1 is extracted then we have the case 9--1(ell).
If the curve 2 is extracted then we have the case 22--1(ell).
If the curve 3 is extracted then we have the case 35(ell).

\par {\bf D).} $\widehat S=S(\AAA_4)$.
Then $r=0$ or 1. If $r=0$ then by lemma \ref{CC1}
we have the case 51--1(ell). Let $r=1$.
Then $\widehat D=b\widehat C+5(1-b)\widehat E_1$. If $b=6/7$ then
$\delta(\widehat S,\widehat D)=2$ (on the diagram the exceptional divisor with
discrepancy $-6/7$ is $(-6)$ curve). Therefore $b>6/7$.
We have the next diagram.
\\
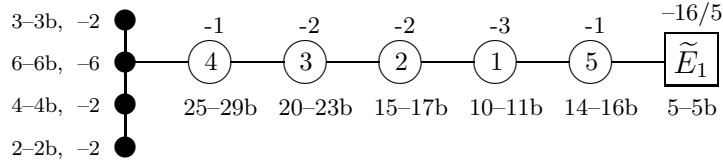
\begin{center}
\begin{picture}(250,70)(0,0)
\put(43,16){\circle*{8}}
\put(25,13){\tiny{--2}}
\put(0,13){\tiny{2--2b,}}
\put(43,20){\line(0,1){8}}
\put(43,32){\circle*{8}}
\put(25,29){\tiny{--2}}
\put(0,29){\tiny{4--4b,}}
\put(43,36){\line(0,1){8}}
\put(43,48){\circle*{8}}
\put(25,45){\tiny{--6}}
\put(0,45){\tiny{6--6b,}}
\put(43,52){\line(0,1){8}}
\put(43,64){\circle*{8}}
\put(25,61){\tiny{--2}}
\put(0,61){\tiny{3--3b,}}
\put(47,48){\line(1,0){20}}
\put(75,48){\circle{16}}
\put(65,28){\scriptsize{25--29b}}
\put(73,45){\footnotesize{4}}
\put(73,59){\scriptsize{-1}}
\put(83,48){\line(1,0){20}}
\put(111,48){\circle{16}}
\put(101,28){\scriptsize{20--23b}}
\put(109,45){\footnotesize{3}}
\put(109,59){\scriptsize{-2}}
\put(119,48){\line(1,0){20}}
\put(147,48){\circle{16}}
\put(137,28){\scriptsize{15--17b}}
\put(145,45){\footnotesize{2}}
\put(145,59){\scriptsize{-2}}
\put(155,48){\line(1,0){20}}
\put(183,48){\circle{16}}
\put(173,28){\scriptsize{10--11b}}
\put(181,45){\footnotesize{1}}
\put(181,59){\scriptsize{-3}}
\put(191,48){\line(1,0){20}}
\put(219,48){\circle{16}}
\put(209,28){\scriptsize{14--16b}}
\put(217,45){\footnotesize{5}}
\put(217,59){\scriptsize{-1}}
\put(227,48){\line(1,0){20}}
\put(247,44){\fbox{$\widetilde E_1$}}
\put(248,28){\scriptsize{5--5b}}
\put(246,64){\scriptsize{--16/5}}
\end{picture}
\end{center}

If the curve with number $k-1$ is extracted then we have the case
$51-k(\rm ell)$.
The curve
1 can be extracted with discrepancy $0$ or $-1/2$.

\par {\bf E).} $\widehat S=S(2\AAA_1+\AAA_3)$.
Then $r=0,1,2$. If $r=0$ then by lemma \ref{CC1}
we have the case 52--1(ell). Let $r=1$.
Then $\widehat D=b\widehat C+4(1-b)\widehat E_1$.
The extraction of required curve is possible only at two points:
$\AAA_1$ and $\AAA_3$. At the point $\AAA_1$ it is absent. At the point
$\AAA_3$ we have the next diagram.
\\
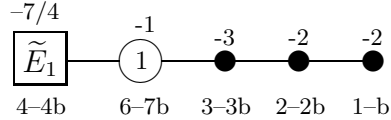
\begin{center}
\begin{picture}(155,40)(0,0)
\put(10,20){\fbox{$\widetilde E_1$}}
\put(11,4){\scriptsize{4--4b}}
\put(9,40){\scriptsize{--7/4}}
\put(30,24){\line(1,0){20}}
\put(58,24){\circle{16}}
\put(50,4){\scriptsize{6--7b}}
\put(56,21){\footnotesize{1}}
\put(56,35){\scriptsize{-1}}
\put(66,24){\line(1,0){20}}
\put(90,24){\circle*{8}}
\put(81,4){\scriptsize{3--3b}}
\put(86,30){\scriptsize{-3}}
\put(94,24){\line(1,0){20}}
\put(118,24){\circle*{8}}
\put(109,4){\scriptsize{2--2b}}
\put(114,30){\scriptsize{-2}}
\put(122,24){\line(1,0){20}}
\put(146,24){\circle*{8}}
\put(138,4){\scriptsize{1--b}}
\put(142,30){\scriptsize{-2}}
\end{picture}
\end{center}

We obtained the case 52--2(ell).
Let $r=2$. Then $a_1+a_2=4(1-b)\le 4/7$, a contradiction
with lemma \ref{vida}.

\par {\bf F).} $\widehat S=S(\DDD_5)$.
Then $r=0$ or 1. If $r=0$ then by lemma \ref{CC1}
we have the case 53--1(ell). Let $r=1$.
Then $\widehat D=b\widehat C+4(1-b)\widehat E_1$. If $b=6/7$ then
$\delta(\widehat S,\widehat D)=2$ (on the diagram the exceptional divisor with
discrepancy $-6/7$ is the central curve). Therefore $b>6/7$.
We have the next diagram.
\\
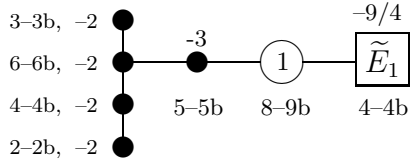
\begin{center}
\begin{picture}(160,70)(0,0)
\put(43,16){\circle*{8}}
\put(25,13){\tiny{--2}}
\put(0,13){\tiny{2--2b,}}
\put(43,20){\line(0,1){8}}
\put(43,32){\circle*{8}}
\put(25,29){\tiny{--2}}
\put(0,29){\tiny{4--4b,}}
\put(43,36){\line(0,1){8}}
\put(43,48){\circle*{8}}
\put(25,45){\tiny{--2}}
\put(0,45){\tiny{6--6b,}}
\put(43,52){\line(0,1){8}}
\put(43,64){\circle*{8}}
\put(25,61){\tiny{--2}}
\put(0,61){\tiny{3--3b,}}
\put(47,48){\line(1,0){20}}
\put(71,48){\circle*{8}}
\put(67,54){\scriptsize{-3}}
\put(62,28){\scriptsize{5--5b}}
\put(75,48){\line(1,0){20}}
\put(103,48){\circle{16}}
\put(101,45){\footnotesize{1}}
\put(95,28){\scriptsize{8--9b}}
\put(111,48){\line(1,0){20}}
\put(131,44){\fbox{$\widetilde E_1$}}
\put(132,28){\scriptsize{4--4b}}
\put(130,64){\scriptsize{--9/4}}
\end{picture}
\end{center}

We obtained the case 53--2(ell).

\par {\bf G).} $\widehat S=S(3\AAA_2)$.
Then $r=0$ or 1. If $r=0$ then by the proof of corollary \ref{corell1} it follows
that $\widetilde C^2\ge 4$. Therefore the variant $r=0$ is impossible.
Let $r=1$. Then there are no cases since there is no an extracted curve
with required discrepancy for the pair
$(\CC^2,3(1-b)\{x=0\})/\ZZ_3(2,1)$.

\par {\bf H).} $\widehat S=S(\AAA_1+\AAA_5)$.
Similarly to point {\bf G)} we have $r=1,2$. Just as in point
{\bf E)} $r\ne 2$.
Let $r=1$.
Then $\widehat D=b\widehat C+3(1-b)\widehat E_1=\frac67\widehat C+
\frac37\widehat E_1$ (see lemma \ref{vida}), where $\widehat E_1$ is one of two
$(-1)$ curves (on a minimal resolution).
The extraction can happen only if
$\widehat E_1$ is $(-1)$ curve not passing through $\AAA_1$.
\\
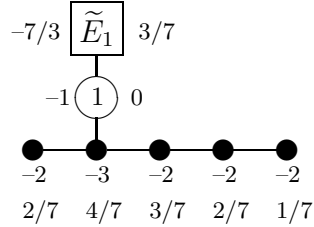
\begin{center}
\begin{picture}(110,84)(0,0)
\put(10,30){\circle*{8}}
\put(6,18){\scriptsize{--2}}
\put(6,4){\scriptsize{2/7}}
\put(14,30){\line(1,0){16}}
\put(34,30){\circle*{8}}
\put(30,18){\scriptsize{--3}}
\put(30,4){\scriptsize{4/7}}
\put(38,30){\line(1,0){16}}
\put(58,30){\circle*{8}}
\put(54,18){\scriptsize{--2}}
\put(54,4){\scriptsize{3/7}}
\put(62,30){\line(1,0){16}}
\put(82,30){\circle*{8}}
\put(78,18){\scriptsize{--2}}
\put(78,4){\scriptsize{2/7}}
\put(86,30){\line(1,0){16}}
\put(106,30){\circle*{8}}
\put(102,18){\scriptsize{--2}}
\put(102,4){\scriptsize{1/7}}
\put(34,34){\line(0,1){8}}
\put(34,50){\circle{16}}
\put(32,47){\footnotesize{1}}
\put(15,47){\scriptsize{--1}}
\put(47,47){\scriptsize{0}}
\put(34,58){\line(0,1){8}}
\put(24,71){\fbox{$\widetilde E_1$}}
\put(2,72){\scriptsize{--7/3}}
\put(50,72){\scriptsize{3/7}}
\end{picture}
\end{center}

We obtained the case 54(ell).

\par {\bf I).} $\widehat S=S(\EEE_6)$.
Similarly to the previous point $r=1$ and $b=6/7$.
Then $\widehat D=\frac67\widehat C+\frac37\widehat E_1$ and
$\delta(\widehat S,\widehat D)=2$. Indeed, after a blow-up of $\CC^3$
with weights $(3,2,2)$ for the pair

$$
\Big(\widehat S,(3/7)\widehat E_1\Big)\an \Big(x^2+y^3+z^3=0\subset \CC^3,
(3/7)\{z=0\}\Big)
$$
the exceptional divisor with discrepancy
$a(\ \cdot\ ,\frac37\widehat E_1)=-6/7$ appears. Therefore there are no cases.

\subsection{} {\it Consider the second case $\rho(\widehat S)=2$ and $r=2$.}
\par
According to proposition \ref{CC} $\widetilde C^2\ge 7$, except the case $n_1=n_2=2$,
$b=6/7$ and $B=0$.
\begin{lemma}\label{help1}
Let $\widetilde C^2\ge 7$. Then the surface
$\widehat S$ can be obtained by usual blow-up of point on the cone $\PP(1,1,2)$, which
is not its vertex.
\begin{proof} Let $\psi\colon \widehat S^{\rm min}\to \widehat S$ be a minimal
resolution and $\overline{S}$ be a minimal model of $\widehat S^{\rm min}$.
All variants of
$\overline{S}$ were described in the proof of proposition \ref{model}. Therefore, if
$\rho(\widehat S^{\rm min})\ge 4$ then we have a contradiction with
$\overline{C}^2\ge 11-\rho(\overline{S})$, where $\overline{C}$ is an
image of $\widehat C$.
\par
If $\rho(\widehat S^{\rm min})=3$ then we have the case in the lemma statement.
Moreover
$\widetilde C^2=7$ and the singular points $\AAA_1$ and $\AAA_2$ lie on $C$.
\par
If $\rho(\widehat S^{\rm min})=2$ then $\widehat S=\FFF_0$ since $\widehat C$ is
very ample divisor. The divisor
$-(K_{\FFF_0}+b\widehat C+(3/7)\widehat E_1+(3/7)\widehat E_2)$ is not nef, a
contradiction.
\end{proof}
\end{lemma}

\begin{corollary} If $\widetilde C^2\ge 7$ then $(S,D)$ is the log surface
$6-2(\rm ell)$.
\end{corollary}

Now let $n_1=n_2=2$, $b=6/7$ and $B=0$. Let us prove that there are no cases.
Similarly to previous lemma
\ref{help1}
$\rho(\widehat S^{\rm min})=2,3,4$ and the variant
$\rho(\widehat S^{\rm min})=2$ is impossible.
Besides, if $\rho(\widehat S^{\rm min})=3$ then
$\overline{C}^2=9-\rho(\overline{S})$. A contradiction.
\par
Let $\rho(\widehat S^{\rm min})=4$. Then $\FFF_2$ is a minimal model of
$\widehat S^{\rm min}$ and
$\varphi \colon \widehat S^{\rm min}\to \FFF_2=\overline{S}.$
There are two possibilities: the proper transforms of $\widehat E_1$ and
$\widehat E_2$ don't contract by
$\varphi$; the proper transform of $\widehat E_1$ is contracted by $\varphi$ but
the proper transform of $\widehat E_2$ is not.
In the first possibility
$-(K_{\FFF_2}+\frac67(2E_0)+(3/7)\overline{E}_1+(3/7)\overline{E}_2)$ is not a nef
divisor. We get a contradiction.
In the second possibility there is an exceptional curve $\Gamma$ on
$\widehat S^{\rm min}$ with $a(\Gamma)=a(\Gamma,(6/7)\widehat C+(3/7)\widehat E_1+
(3/7)\widehat E_2)<0$.

\begin{lemma} $a(\Gamma)\le -3/14$.
\begin{proof} If $\psi(\Gamma)\in \Supp \widehat D$ then in the same way as in
lemma \ref{vida} we obtain the required inequality.
If $\psi(\Gamma)\notin \Supp \widehat D$ then
$\psi(\Gamma)\in \widehat S$ is not Du Val singularity.
Therefore $a(\Gamma)\le -1/3$ \cite[proposition 2.4.8]{KK}.
\end{proof}
\end{lemma}

The divisor
$-(K_{\FFF_2}+\frac67(2E_0)+(3/7)\overline{E}_2+(3/14)\overline{\Gamma})$
is not nef. We get a contradiction.

\section{\bf Case ${\bf p_a(C)=0} $ and ${\bf \widetilde C^2\ge 1}$}

Let $p_a(C)=0$.
The next proposition is proved as previous
proposition \ref{CC}.
\begin{proposition}\label{CC2}
$\deg\Diff_C(0)=\sum_{i=1}^r (1-1/n_i)$ and
$$C^2=\frac{-2+\deg \Diff_C(0)+\sum b_i(B_i\cdot C)}{1-b},
\ \ \ \widetilde C^2=C^2-\sum_{i=1}^r\frac{q_i}{n_i}.$$
\end{proposition}

\begin{corollary} $\widetilde C^2\ge -1$.
\begin{proof} By proposition \ref{CC2} $\widetilde C^2\ge 7(-2+\Diff_C(0)+
\sum b_i(B_i\cdot C))-\Diff_C(0)=
6(-2+\Diff_C(0)+\sum b_i(B_i\cdot C))-2+\sum b_i(B_i\cdot C)
>-2+\sum b_i(B_i\cdot C)\ge -2$.
\end{proof}
\end{corollary}

\par
Consider the first case
$\widetilde C^2\ge 1$.
By proposition \ref{free2} a linear system $|\widetilde C|$ gives a birational
morphism $h$ and
$\widehat S=\PP(1,1,n)$, $\FFF_n$.

\begin{theorem} Let $(\widehat S,\widehat C)\simeq(\FFF_n,E_0+kf)$. Then
$n=1,2,3$; $k=1$  and we obtain four cases
$2-2(+4)$, $3-1(+5)$, $6-3(+3)$, $8-2(+3)$.
\begin{proof} Let $n\ne 0$. Since $\widehat E_i\cdot \widehat C=1$ then either
$\widehat E_i\sim f$ for all $i$, or
$\widehat E_1= E_{\infty}$, $\widehat E_i\sim f$ for all $i\ge 2$ and $k=1$.
\par
Consider the first variant.
Then there exists 1-, 2-, 3-, 4- or 6-complement
$D^+\ge \widehat C+\sum b_i\widehat B_i$ of the divisor
\begin{equation}
K_{\widehat S}+b\widehat C+\sum_{i=1}^ra_i\widehat E_i+\sum b_i\widehat B_i
\sim_{\QQ} 0
\end{equation}
in the neighborhood of fiber $f'$ by theorem \ref{contr}.
From (1) we have $(K_{\widehat S}+D^+)\cdot f'>0$, a contradiction.
\par
Consider the second variant. Then the divisor from (1) is $\QQ$-linear equivalent
to the next one
\begin{equation*}
K_{\widehat S}+b(E_0+f)+a_1E_{\infty}+ \sum_{i\ge 2}a_if+\sum b_i(k_iE_0+l_if)
\sim_{\QQ} 0.
\end{equation*}
It is obvious that there exists a number $i_0$ such that $k_{i_0}\ne 0$. Since
$0\le -(K_{\widehat S}+(6/7)(E_0+f)+(1/2)E_0)\cdot E_0=8/7-(5/14)n$ then $n\le 3$.
\par
Let $n=3$. Then $\widehat D=(7/8)\widehat C+(5/8)\widehat E_1+(1/2)B_1$, where
$B_1\sim E_0$. Thus $r=1$, $\widetilde S=\widehat S$. We obtain the case
$3-1(+5)$.

\par
Let $n=2$. Then $\widehat D=(8/9)\widehat C+(4/9)\widehat E_1+(2/3)B_1$ or
$\widehat D=(6/7)\widehat C+(9/14)\widehat E_1+(3/7)\widehat E_2+(1/2)B_1$,
where $B_1\sim E_0$. In the first variant $r=1$ and we obtain the case $2-2(+4)$.
In the second variant $r=2$.
The extraction of required curve can take place at the point
$\widehat E_1\cap \widehat E_2$ only. It is easy to prove that it is absent.

\par
Let $n=1$. Then $r\le 2$ and
$\widehat D=(6/7+a)\widehat C+(9/14-a)\widehat E_1+(11/14-2a)\widehat E_2+(1/2)B_1$
or
$\widehat D=(6/7+a)\widehat C+(10/21-a)\widehat E_1+(13/21-2a)\widehat E_2+(2/3)B_1$,
where $B_1\sim E_0$ and $a\in[0,\frac17)$. In both variants $r=2$
and the extraction of required curve can take place at the point
$\widehat E_1\cap \widehat E_2$ only. For the pair
$(\CC^2_{x,y},(9/14-a)\{x=0\}+(11/14-2a)\{y=0\})$ we can extract a divisor
with discrepancy 0 only.
It is not hard to prove that it happens by the weighted blow-ups with weights
$(\alpha,\beta)=(1,2),(1,3)$ and $(2,1)$. We obtain the cases
$6-3(+3)$ and $8-2(+3)$ in the first and second possibilities for
$(\alpha,\beta)$ respectively.
The possibility $(2,1)$ is not realized since
a proper transform of
$\widehat E_2$ is $(-1)$ curve on a minimal resolution.
For the same reason ($\beta=1$) the second variant of
$\widehat D$ is not realized.
\par
Let $n=0$. Then the number of intersection points between
$\widehat E_i$ is not more then
$r-1$.
Therefore, one can assume without loss of generality that
$\widehat E_1\sim E_0$, $\widehat E_i\sim f$ for $i\ge 2$.
By the same argument we can prove that it is impossible to convert every
$\widehat E_i$ to $(-k_i)$ curve (on a minimal resolution), where
$k_i\ge 2$. Hence, there are no cases.
\end{proof}
\end{theorem}

\begin{proposition}
Let $(\widehat S,\widehat C)=(\PP(1,1,n),X_n)$. Then $n\le 5$ and
$\widehat E_i$ are the generators of cone.
\begin{proof}
Since $\widehat E_i\cdot \widehat C=1$ then
$\widehat E_i$ is the generator of cone.
It is clear that  $K_{\widehat S}+\widehat D$ is
$\frac17$-log terminal divisor at the cone vertex.
Let $n\ge 6$ and
$\FFF_n\to \widehat S$ be a minimal resolution. A proper transform of
$\widehat B_i$ is denoted by $\overline{B}_i$. Then
$\overline{B}_i\sim f$ for all $i$. Indeed, if
$\overline{B}_i\sim l_iE_0+k_if$, where $l_i\ge 1$ then
$-(K_{\widehat S}+(6/7)\widehat C+(1/2)\widehat B_i)$ is not a nef divisor,
a contradiction. Hence
$K_{\widehat S}+\widehat D$ is not $\frac17$-log terminal divisor at the cone vertex.
\end{proof}
\end{proposition}

Consider case by case all variants. The requirement to be $\frac17$-log terminal at
the cone vertex ($n\ge 2$) implies the existence of $i_0$ such that
$(\widehat B_{i_0}-X_n)$ is a nef divisor. Put $i_0=1$.

\par {\bf A).} $\widehat S=\PP(1,1,5)$. We have $r=0$ and the case $5(+5)$.

\par {\bf B).} $\widehat S=\PP(1,1,4)$. Then $r=0,1$ and if $r=0$ then
we have the case $4(+4)$. Let $r=1$. Then $\widehat D=b\widehat C+(1/2)\widehat B_1+
(4-4b)\widehat E_1$, where $\widehat B_1=X_4$.
The extraction of required curve can take place at the point
$\widehat E_1\cap \widehat B_1$ and at the cone vertex only. It can easily be checked
that it is absent.

\par {\bf C).} $\widehat S=\PP(1,1,3)$. We have the following possibilities for
$\widehat D$:
\begin{enumerate}
\item $\widehat D=(\frac53-\frac{k-1}k)\widehat C+\frac{k-1}k\widehat B_1$,
where $k=4,5$. We obtain the case $3-2(+3)$.
\item $\widehat D=\frac67\widehat C+\frac23\widehat B_1+\frac37\widehat E_1$.
It is not realized.
\item $\widehat D=\frac67\widehat C+\frac12\widehat B_1+
\frac12\widehat B_2+\frac37\widehat E_1$. We obtain the case $10-1(+3)$.
\item $\widehat D=b\widehat C+\frac12\widehat B_1+
a_1\widehat E_1+a_2\widehat E_2$. It is not realized.
\item $\widehat D=(\frac56+\frac1{3k})\widehat C+\frac12\widehat B_1+
\frac{k-1}k\widehat B_2$, where $k=3,4,5,6$. It is the case $3-3(+3)$.
\item $\widehat D=b\widehat C+\frac12\widehat B_1+
(\frac72-3b)\widehat E_1$. We obtain the cases
$7-1(+3)$, $14-1(+3)$, $10-2(+3)$, $20(+3)$, $29(+3)$.
\end{enumerate}

Here $\widehat B_1=X_3$, $\widehat B_2=X_1$, $a_1+a_2=7/2-3b$.
In the third possibility it is possible that
$\widehat B_3=\widehat B_1+\widehat B_2$ is an irreducible curve $X_4$.
\par
It follows easily that possibilities (2),(4) are not realized.
In possibility (3) the extraction can be only if
$\widehat B_3\cap \widehat E_1=P$, where
$P$ is a cone vertex. It is the case $10-1(+3)$.
The calculations are similarly as in possibility (6).
\par
In possibility (6) $r=1$ and the extraction of required curve can take place at
the point
$Q=\widehat E_1\cap\widehat B_1$ and at the cone vertex. Consider the point $Q$.
Let $h\colon\widetilde S\to (\widehat S\ni Q)$ be an required extraction.
Then condition $\delta(S,D)=1$ gives the following requirement: on a minimal
resolution of surface $\widetilde S$ a proper transform of $\widetilde E_1$ is a curve
with self-intersection index $k=-2,-3$ (see proposition \ref{vidc}) and there is
only one singular point of $\widetilde S$, which is the cone vertex and lies on
$\widetilde E_1$.
For the pair $(\CC^2_{x,y}, \frac12\{x=0\}+(\frac72-3b)\{y=0\})$
these two variants correspond to the weighted blow-ups with weights
(1,2) and (1,3). We have the cases
$7-1(+3)$, $14-1(+3)$.
\par
Consider the cone vertex. Similarly, by proposition
\ref{vidc} there are no singular points of $\widetilde S$ which lie
on $\widetilde E_1$ and
$\widetilde E_1^2=-2,-3,-4$.
The realization of these three possibilities is shown in the next diagram.
\\
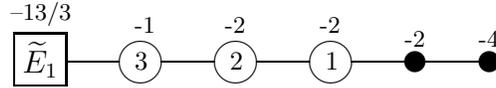
\begin{center}
\begin{picture}(196,40)(0,0)
\put(10,20){\fbox{$\widetilde E_1$}}
\put(9,40){\scriptsize{--13/3}}
\put(30,24){\line(1,0){20}}
\put(58,24){\circle{16}}
\put(56,21){\footnotesize{3}}
\put(56,35){\scriptsize{-1}}
\put(66,24){\line(1,0){20}}
\put(94,24){\circle{16}}
\put(92,21){\footnotesize{2}}
\put(90,35){\scriptsize{-2}}
\put(102,24){\line(1,0){20}}
\put(130,24){\circle{16}}
\put(128,21){\footnotesize{1}}
\put(126,35){\scriptsize{-2}}
\put(138,24){\line(1,0){20}}
\put(162,24){\circle*{8}}
\put(158,30){\scriptsize{-2}}
\put(166,24){\line(1,0){20}}
\put(190,24){\circle*{8}}
\put(186,30){\scriptsize{-4}}
\end{picture}
\end{center}

We obtain the cases $10-2(+3)$, $20(+3)$, $29(+3)$.

\par {\bf D).} $\widehat S=\PP(1,1,2)$.
We have the following possibilities for $\widehat D$. Here
$\widehat B_1=X_2$, $\widehat B_2=X_1$, $\widehat B_3=X_3$ and $\widehat B_4=X_1$.
\par (1). $\widehat D=\frac{11}{12}\widehat C+\frac56\widehat B_1+
\frac12\widehat B_2$. It is the case $2-3(+2)$.
\par (2). $\widehat D=b\widehat C+\frac56\widehat B_1+
(\frac73-2b)\widehat E_1$. The extraction can take place at the point
$\widehat B_1\cap\widehat E_1$ only.
Now let us give a reasoning which will be used many times. It allows to decrease
the computation quantity.
Thus, consider the pair $(\CC^2_{x,y},(\frac73-2b)\{x=0\}+\frac56\{y=0\})$.

\begin{lemma}\label{weight} The extraction of required curve can take place
by a blow-up with weights
$(2,1)$ only. It is the case $6-4(+2)$.
\begin{proof} It can easily be checked that the extraction of required curve
with a discrepancy $1/\vartheta-1$ can be realized by a toric blow-up only, i.e.
by a weighted blow-up with weights
$(\alpha,\beta)$. Take a blow-up $\widetilde S\to S$ with weights $(\alpha,\beta)$.
The condition that $\widetilde E_1$ is not  $(-1)$ curve on a minimal
resolution
$\widetilde S$ is $\widetilde E_1^2<-1$, i.e.
$-\alpha/\beta<-1$. Hence $\alpha\ge \beta+1$. Solving the equation system

$$
\left\{
\begin{array}{l}
\alpha(12b-8)+\beta=\frac6{\vartheta}  \\
\alpha\ge \beta+1\\
\alpha, \beta, \vartheta \in \ZZ_{>0};\ (\alpha, \beta)=1;\ b\in[\frac67,1),\\
\end{array}
\right.
$$
We get $(\alpha,\beta)=(2,1)$, $\vartheta=1$, $b=7/8$.
\end{proof}
\end{lemma}

\par (3). $\widehat D=(\frac65-\frac{k-1}{2k})\widehat C+\frac45\widehat B_1+
\frac{k-1}k\widehat B_2$, where $k=2,3$. It is the case $2-4(+2)$.

\par (4). $\widehat D=b\widehat C+\frac45\widehat B_1+
(\frac{12}5-2b)\widehat E_1$.
The extraction can happen at the point
$\widehat B_1\cap\widehat E_1$ only. Likewise (2) the extraction is realized by
a blow-up with weights
(2,1) and $b=9/10$, $\vartheta=1$. We obtain the case $6-4(+2)$.

\par (5). $\widehat D=(\frac54-\frac{k-1}{2k})\widehat C+\frac34\widehat B_1+
\frac{k-1}k\widehat B_2$, where $k=3,4$. It is the case $2-5(+2)$.

\par (6). $\widehat D=b\widehat C+\frac34\widehat B_1+
(\frac52-2b)\widehat E_1$.
The extraction can happen at the point
$\widehat B_1\cap\widehat E_1$ only. Likewise (2) the extraction is realized by
the blow-ups with weights
(2,1), (3,1)
and $\vartheta=1$. We obtain the cases $6-4(+2)$ and $9-2(+2)$.

\par (7). $\widehat D=(\frac43-\frac{k-1}{2k})\widehat C+\frac23\widehat B_1+
\frac{k-1}k\widehat B_2$, where $k=4,5,6$. It is the case $2-6(+2)$.

\par (8). $\widehat D=b\widehat C+\frac23\widehat B_1+
(\frac83-2b)\widehat E_1$. Consider the extraction at the cone vertex.
The condition $\delta(S,D)=1$ implies (see proposition \ref{vidc})
$\widetilde E_1\cap \Sing\widetilde S=\emptyset$ and
$\widetilde E_1^2=-2,-3$. These cases are realized.
\\
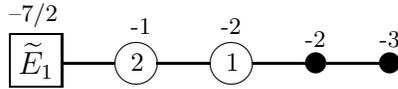
\begin{center}
\begin{picture}(150,40)(0,0)
\put(10,20){\fbox{$\widetilde E_1$}}
\put(10,40){\scriptsize{--7/2}}
\put(30,24){\line(1,0){20}}
\put(58,24){\circle{16}}
\put(56,21){\footnotesize{2}}
\put(56,35){\scriptsize{-1}}
\put(66,24){\line(1,0){20}}
\put(94,24){\circle{16}}
\put(92,21){\footnotesize{1}}
\put(90,35){\scriptsize{-2}}
\put(102,24){\line(1,0){20}}
\put(126,24){\circle*{8}}
\put(122,30){\scriptsize{-2}}
\put(130,24){\line(1,0){20}}
\put(154,24){\circle*{8}}
\put(150,30){\scriptsize{-3}}
\end{picture}
\end{center}

We obtain the cases $7-2(+2)$, $13-1(+2)$.
In the last case, if
$t=19/21$ then $\delta(S,D)=2$.
\par

Consider the extraction at the point
$\widehat B_1\cap\widehat E_1$.
The condition $\delta(S,D)=1$ implies (see proposition \ref{vidc})
$\widetilde E_1\cap \Sing\widetilde S=P$, where $P$ is a cone vertex, $\vartheta=1$ and
$\widetilde E_1^2=-2+\frac12,-3+\frac12,-4+\frac12$.
The case $\widetilde E_1^2=-2+\frac12$ is not realized since
$b<1$. Other two cases are realized by the blow-ups with weights
$(3,1)$ and $(4,1)$.
We obtain the cases $9-2(+2)$, $15-1(+2)$.

\par (9). $\widehat D=b\widehat C+\frac23\widehat B_1+
a_1\widehat E_1+a_2\widehat E_2$, where $a_1+a_2=8/3-2b$.
It is not hard to prove that this possibility is not realized.

\par (10). $\widehat D=(\frac54-\frac{k-1}{2k})\widehat C+\frac12\widehat B_3+
\frac{k-1}k\widehat B_2$, where $k=3,4$. It can happen that
$\widehat B_3=\widehat B_1+\widehat B_4$.
It is the case $2-7(+2)$.

\par (11). $\widehat D=b\widehat C+\frac12\widehat B_3+
(\frac52-2b)\widehat E_1$. It is possible that
$\widehat B_3=\widehat B_1+\widehat B_2$.
\par
Let $\widehat B_3\cap\widehat E_1\ne P$, where $P$ is a cone vertex. If
$\widehat B_3=\widehat B_1+\widehat B_2$ then this condition is always fulfilled.
Considering the cone vertex as in point (8) we obtain the cases
$7-3(+2)$, $13-2(+2)$.
Consider the smooth point $\widehat S$, where $\widehat B_3$ intersects
$\widehat E_1$. Similarly the extraction can take place by a blow-up with
weights (2,1) only. It is the case $6-5(+2)$.
\par
Let $\widehat B_3\cap\widehat E_1= P$. Then
\\
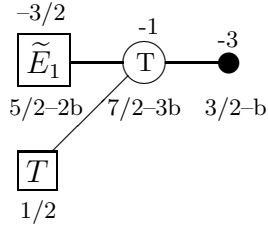
\begin{center}
\begin{picture}(150,80)(0,-40)
\put(10,20){\fbox{$\widetilde E_1$}}
\put(7,4){\scriptsize{5/2--2b}}
\put(10,40){\scriptsize{--3/2}}
\put(30,24){\line(1,0){20}}
\put(58,24){\circle{16}}
\put(44,4){\scriptsize{7/2--3b}}
\put(55,21){\footnotesize{T}}
\put(56,35){\scriptsize{-1}}
\put(66,24){\line(1,0){20}}
\put(90,24){\circle*{8}}
\put(81,4){\scriptsize{3/2--b}}
\put(86,30){\scriptsize{-3}}
\put(10,-21){\fbox{$T$}}
\put(11,-35){\scriptsize{1/2}}
\put(23,-10){\line(1,1){29}}
\end{picture}
\end{center}

Later on the required curve is extracted by a weighted blow-up at the point
$\widetilde E_1\cap T$. Let its discrepancy is equal to
$1/\vartheta-1$.
The condition $\delta(S,D)=1$ implies $7/2-3b<6/7$, i.e. $b>37/42$. Similarly to lemma
\ref{weight} we have to solve the equation $\alpha(4b-3)+\beta(6b-5)=2/\vartheta$.
We obtain $\vartheta=2$ and $(\alpha,\beta)=(1,1)$; $\vartheta=1$ and
(1,2), (1,3), (1,4),
(1,5), (2,1), (2,3), (3,1). We have the cases $7-4(+2)$, $14-2(+2)$, $26(+2)$,
$36(+2)$, $40(+2)$, $13-3(+2)$, $34(+2)$,
$21-1(+2)$ respectively.

\par (12). $\widehat D=b\widehat C+\frac12\widehat B_1+ \frac{k-1}k\widehat B_2+
(2+\frac1k-2b)\widehat E_1$, where $k=3,4,5,6$.
It is not hard to prove that this possibility is not realized.

\par (13). $\widehat D=b\widehat C+\frac12\widehat B_1+
a_1\widehat E_1+a_2\widehat E_2$, where $a_1+a_2=3-2b$.
It is not hard to prove that this possibility is not realized.

\par (14). $\widehat D=\frac67\widehat C+\frac12\widehat B_1+
\frac37\widehat E_1+\frac37\widehat E_2+\frac37\widehat E_3$.
It is not hard to prove that this possibility is not realized.

\par {\bf E).} $(\widehat S,\widehat C)=(\PP^2,\OO_{\PP^2}(2))$. Since
$\widehat E_i\cdot \widehat C=1$ then $r=0$. It is the case $1-1(+4)$.

\par {\bf F).} $(\widehat S,\widehat C)=(\PP^2,\OO_{\PP^2}(1))$.
Consider case by case all possibilities for $\widehat D$.

\par (I). Let $\Supp \widehat D$ consists of the straight lines only and
$s$ is their quantity.
If $r=0$ then we have the cases $1-2(+1)$, $1-3(+1)$, $1-4(+1)$.
\par
Let $r=1$. Then $s=4,5$. At first consider $s=5$.
Then by proposition \ref{vidc} three straight lines pass through
some point on
$\widehat E_1$, which doesn't lie on $\widehat C$. Let
$\widehat B_1$ and $\widehat B_2$ be among them.
Since $b_1+b_2+b_3=3-b-a_1\le 12/7$ (see lemma \ref{vida}) then
$(b_1,b_2,b_3)=(\frac12,\frac12,\frac{k-1}k)$, where $k=2,3$ or
$(b_1,b_2,b_3)=(\frac12,\frac23,\frac12)$. Since
$\widetilde E_1^2\le -2$ then $b_1+b_2+3a_1-3\ge 0$.
Substituting $b_1$, $b_2$, $b_3$ we obtain a contradiction. Now let $s=4$.
Then all straight lines from $\Supp \widehat D$ are in the general position.
By proposition \ref{vidc}
$\widetilde E_1^2=-2,-3,-4$ and $\widetilde E_1\cap \Sing \widetilde S=\emptyset$, i.e.
we have the extractions by the blow-ups with weights
$(\alpha,1)$ only, where $\alpha=3,4,5$.
We obtain the cases $6-6(+1)$, $8-3(+1)$, $11-1(+1)$.

\par
Let $r=2$. Then $s=4,5$. For the same reason the case $s=5$ is impossible.
Thus $s=4$ and all straight lines from $\Supp \widehat D$ are in the general position.
We have two extraction types, where
$\widetilde\Gamma_i$ are the extracted curves.
\\
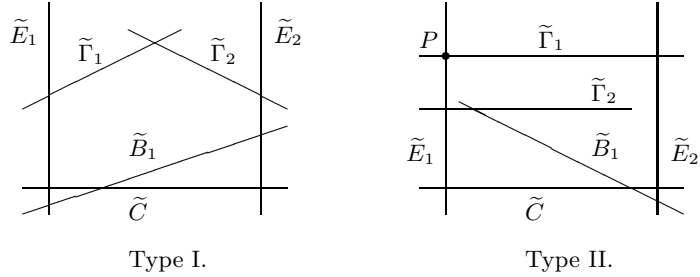
\begin{center}
\begin{picture}(260,100)(0,0)
\put(10,30){\line(1,0){100}}
\put(50,0){\scriptsize{Type I.}}
\put(50,18){\scriptsize{$\widetilde C$}}
\put(20,20){\line(0,1){80}}
\put(100,20){\line(0,1){80}}
\put(5,85){\scriptsize{$\widetilde E_1$}}
\put(105,85){\scriptsize{$\widetilde E_2$}}
\put(10,20){\line(3,1){100}}
\put(50,42){\scriptsize{$\widetilde B_1$}}
\put(10,60){\line(2,1){60}}
\put(110,60){\line(-2,1){60}}
\put(31,78){\scriptsize{$\widetilde\Gamma_1$}}
\put(81,78){\scriptsize{$\widetilde\Gamma_2$}}
\put(160,30){\line(1,0){100}}
\put(200,0){\scriptsize{Type II.}}
\put(200,18){\scriptsize{$\widetilde C$}}
\put(170,20){\line(0,1){80}}
\put(250,20){\line(0,1){80}}
\put(155,40){\scriptsize{$\widetilde E_1$}}
\put(255,40){\scriptsize{$\widetilde E_2$}}
\put(160,80){\line(1,0){100}}
\put(170,80){\circle*{3}}
\put(160,83){\scriptsize{$P$}}
\put(205,83){\scriptsize{$\widetilde\Gamma_1$}}
\put(160,60){\line(1,0){80}}
\put(225,63){\scriptsize{$\widetilde\Gamma_2$}}
\put(260,20){\line(-2,1){85}}
\put(225,42){\scriptsize{$\widetilde B_1$}}
\end{picture}
\end{center}

\par
Consider the type I.
By proposition \ref{vidc}
$\widetilde E_i^2=-n_i$, where $n_i=2,3,4$,
$\widetilde E_i\cap \Sing \widetilde S=\emptyset$ and the discrepancies
$a(\widetilde\Gamma_i,\widehat D)=0$.
Thus
\begin{equation}
3=b+b_1+a_1+a_2=b+b_1+1+\frac{b+b_1-2}{n_1}+1+\frac{b+b_1-2}{n_2}.
\end{equation}
Taking into account $6/7\le b<1$ we have $b_1=1/2$ and $(n_1,n_2)=(2,3)$, (2,4),
(3,3), (3,4).
Two required toric blow-ups have the weights
$(n_1+1,1)$ and $(1,n_2+1)$.
We obtain the cases $27(+1)$, $43(+1)$, $44(+1)$, $42(+1)$ respectively.

\par
Consider the type II.
By proposition \ref{vidc}
$\widetilde E_2^2=-(m-1)$, where $m=3,4,5$;
$\widetilde E_2\cap \Sing \widetilde S=\emptyset$;
$\widetilde E_1\cap \Sing \widetilde S=P$
and the discrepancies $a(\widetilde\Gamma_i,\widehat D)=0$.
Thus, at the point $P$ the surface $\widetilde S$ has a singularity $\AAA_{m-1}$
and the curve
$\widetilde E_1$ is contracted to a singular point $\frac1{km-m+1}(m,1)$, where
$k=-\widetilde E_1^2+\frac{m-1}m$. The condition $\delta(S,D)=1$ implies $m(k-2)<6$.
Hence
$k=2,3$ and $m=3,4,5$. Taking into account the condition as
(2) we obtain
$17-1(+1)$ and (3,3,1/2), $12-1(+1)$ and (4,2,1/2), $45(+1)$ and (4,3,1/2),
$46(+1)$ and (5,2,1/2), $47-1(+1)$ and (3,2,2/3), where the triples of numbers
are equal to $(m,k,b_1)$.

\par
Let $r=3$. Then $s=4,5$. For the same reason the case $s=5$ is impossible.
Thus $s=4$ and all straight lines from $\Supp \widehat D$ are in the general position.
The extraction of three curves is shown in the next figure.
\\
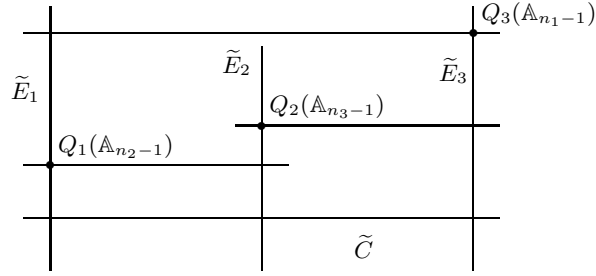
\begin{center}
\begin{picture}(210,110)(0,0)
\put(10,30){\line(1,0){180}}
\put(20,10){\line(0,1){100}}
\put(5,75){\scriptsize{$\widetilde E_1$}}
\put(100,10){\line(0,1){85}}
\put(180,10){\line(0,1){100}}
\put(10,50){\line(1,0){100}}
\put(20,50){\circle*{3}}
\put(23,55){\scriptsize{$Q_1(\AAA_{n_2-1})$}}
\put(90,65){\line(1,0){100}}
\put(100,65){\circle*{3}}
\put(103,70){\scriptsize{$Q_2(\AAA_{n_3-1})$}}
\put(10,100){\line(1,0){180}}
\put(180,100){\circle*{3}}
\put(183,105){\scriptsize{$Q_3(\AAA_{n_1-1})$}}
\put(167,82){\scriptsize{$\widetilde E_3$}}
\put(85,85){\scriptsize{$\widetilde E_2$}}
\put(135,15){\scriptsize{$\widetilde C$}}
\end{picture}
\end{center}

By proposition \ref{vidc} all curves are extracted with discrepancy 0.
Thus $\Sing \widetilde S=Q_1\cup Q_2\cup Q_3$ and the self-intersection
indexes of proper transforms of $\widetilde E_i$ are equal to $(-n_i)$
on a minimal resolution.
At the points $Q_i$ we have the singularities $\AAA_{n_i-1}$ respectively.
Put $n_1\ge n_2\ge n_3$.
The condition $\delta(S,D)=1$ implies
$$
\left\{
\begin{array}{l}
n_1(n_2-2)<6  \\
n_2(n_3-2)<6  \\
n_3(n_1-2)<6.  \\
\end{array}
\right.
$$

Taking into account the condition as (2) we obtain $(n_1,n_2,n_3)=(3,2,2)$ or
$(4,2,2)$.
We have the cases $18-1(+1)$, $25-1(+1)$.

\par (II). $\Supp \widehat D$ contains an irreducible conic.
If $r=0$ then we have the case $1-4(+1)$.
\par
Let $r=1$. Then we have the following possibilities for $\widehat D$, where
$\widehat B_1=X_2$ and $\widehat B_2=X_1$.

\par (1). $\widehat D=b\widehat C+\frac12\widehat B_1+\frac12\widehat B_2+
(\frac32-b)\widehat E_1$.
It can easily be checked that the extraction can be only if
$\widehat E_1$ is tangent $\widehat B_1$ and
$\widehat B_2$ passes through the point of tangency. We have the pair
$(\CC^2,(\frac32-b)\{x=0\}+\frac12\{x+y^2=0\}+\frac12\{y=0\})$.
Similarly to lemma \ref{weight} we obtain $(\alpha,\beta)=(4,1)$, (5,2), (7,3)
and $\vartheta=1$. These are the cases $8-4(+1)$, $9-3(+1)$, $15-2(+1)$.

\par (2). $\widehat D=b\widehat C+\frac12\widehat B_1+\frac23\widehat B_2+
(\frac43-b)\widehat E_1$. Similarly
$\widehat E_1$ is tangent $\widehat B_1$ and
$\widehat B_2$ passes through the point of tangency. We obtain $(\alpha,\beta)=(3,1)$
and
$\vartheta=1$, i.e. it is the case $6-7(+1)$.

\par (3). $\widehat D=b\widehat C+\frac23\widehat B_1+(\frac53-b)\widehat E_1$.
Let $\widehat E_1$ is not tangent $\widehat B_1$.
Then taking a blow-up of their intersection point we get
$(\alpha,\beta)=(3,1)$ and
$\vartheta=1$, i.e. it is the case $6-8(+1)$.
\par
Let $\widehat E_1$ is tangent $\widehat B_1$. Take a blow-up with weights (2,1). Then
a discrepancy is equal to $a(\ \cdot\ ,\widehat D)=2b-8/3>-6/7$. Therefore $b>19/21$.
Taking it into account we obtain
(3,1) and $\vartheta=2$; (5,1) and $\vartheta=1$; (7,3) and $\vartheta=1$;
(8,3) and $\vartheta=1$; (9,4) and $\vartheta=1$; (11,5) and $\vartheta=1$ as in
lemma \ref{weight}. These are the cases
$6-9(+1)$, $11-2(+1)$, $15-3(+1)$, $22-2(+1)$, $24(+1)$, $30(+1)$
respectively.

\par (4). $\widehat D=b\widehat C+\frac34\widehat B_1+(\frac32-b)\widehat E_1$.
Similarly $\widehat E_1$ is tangent $\widehat B_1$. Then $\vartheta=1$ and
$(\alpha,\beta)=(4,1)$, (5,2), (7,3).
These are the cases $8-5(+1)$, $9-4(+1)$, $15-3(+1)$ respectively.

\par (5). $\widehat D=b\widehat C+\frac{k-1}k\widehat B_1+(1+\frac2k-b)\widehat E_1$,
where $k=5,6$.
Similarly $\widehat E_1$ is tangent $\widehat B_1$. Then $\vartheta=1$ and
$(\alpha,\beta)=(3,1)$.
It is the case $6-10(+1)$.

\par Let $r=2$. Then
$\widehat D=b\widehat C+\frac12\widehat B_1+a_1\widehat E_1+
a_2\widehat E_2$, where $a_1+a_2=2-b$.
Similarly
$\widehat E_1$ is tangent $\widehat B_1$ and
$\widehat E_2$ passes through the point $P$ of tangency.
Besides, two extracted curves are contracted to the point
$P$. Similarly to lemma
\ref{weight} we obtain that these blow-ups correspond to ones with weights
(3,1), (1,3) and the discrepancies are equal to 0.
It is the case $48(+1)$.

\par (III). $\Supp \widehat D$ contains an irreducible cubic $B'$.
Then
$\widehat D=b\widehat C+\frac12B'+(\frac32-b)\widehat E_1$.
By proposition \ref{vidc} $\widehat E_1$ intersects $B'$ at one or two points.
It is easy to prove that the case of two points is not realized.
The case of one point was developed in
II (1), where the cubic was decomposed .
In the obtained answer these new cases were included.
Notice that the condition $\delta(S,D)=1$ implies $b>\frac{37}{42}$ in the case
$(\CC^2,\frac12\{x^2+y^3=0\}+(\frac32-b)\{x=0\})$.
Then $(\vartheta,\alpha,\beta)=(1,5,2)$. It is the case $9-3(+1)$, moreover
the monomial $x_1x_2x_3$ is absent in the polynomial defining $X_9$.

\section{\bf Case ${\bf p_a(C)=0}$ and ${\bf \widetilde C^2=0}$.}

Let $p_a(C)=0$ and $\widetilde C^2=0$.
Let $\psi\colon S'\to S$ be a partial resolution of $S$ along $C$ and
$$K_{S'}+bC'+\sum b_iB_i'+\sum_{i=1}^ra_iE_i'=\psi^*(K_S+D').$$
Let $\widetilde S\to S'$ be a partial resolution of $S'$ along $E_1'$.
Then we have
$$
K_{\widetilde S}+b\widetilde C+\sum b_i\widetilde B_i+
\sum_{i=1}^ra_i\widetilde E_i+\gamma\widetilde\Gamma=f^*(K_S+D').
$$

If $\Sing S'\cap E_1'=\emptyset$ then suppose that $S'=\widetilde S$.
For $m\gg 0$ a linear system $|\widetilde E_1+m\widetilde C|$ gives a birational
morphism $h\colon \widetilde S\to \widehat S$ by
proposition \ref{free2} .
Note that $h$ doesn't contract
$\widetilde E_i$, $\widetilde C$, $\widetilde \Gamma$. We have
$$K_{\widetilde S}+b\widetilde C+\sum b_i\widetilde B_i+
\sum_{i=1}^ra_i\widetilde E_i+\gamma\widetilde\Gamma=h^*(K_{\widehat S}+
b\widehat C+\sum b_i\widehat B_i+
\sum_{i=1}^ra_i\widehat E_i+\gamma\widehat\Gamma).$$
It is clear that $\widehat S=\FFF_n$, where $n\ge 2$;
$\widehat C\sim\widehat \Gamma\sim f$;
$\widehat E_1= E_{\infty}$; $\widehat E_i\sim E_0$ for $i\ge 2$.

\begin{lemma}\label{vida2}
If $S'\ne \widetilde S$ then $a_1\ge 4/7$ and $\gamma\ge 2/7$.
\begin{proof} The proof is the same one as in lemma \ref{vida}.
Note that the equality takes place in the case
$\AAA_2$ and $b=6/7$ only.
\end{proof}
\end{lemma}

\begin{proposition} $1\le r \le 3$.
\begin{proof} Let $r\ge 4$. Then by proposition \ref{CC2} we have $r=4$ and
the singular point $\AAA_1$ lies on
$C$. Let it be the first point $P_1$.
Then $S'=\widetilde S$, $\rho(\widetilde S)=5$ and $\widehat S=\FFF_2$. By proposition
\ref{CC2} there is $j$ such that $a_j\ge 4/7$. Hence $r=4$ and
$\widehat D=b\widehat C+\sum_{i=1}^4a_i\widehat E_i$. The set
$\Phi=\{E_i\cap E_j\mid i,j\ge 1, i>j\}$ consists of at most three different points.
There are four variants $(2P_1,2P_2,2P_3)$, $(3P_1,3P_2)$, $(4P_1,P_2,P_3)$,
$(6P_1)$, where $P_i\in \Phi$ and the numbers at $P_i$ are equal to
$\sum_{j:j\ne i}(E_i\cdot E_j)_{P_i}$.
Considering these variants as in lemma
\ref{weight} we get a contradiction, for instance with $a_2+a_3+a_4=2-a_1\le 11/7$.
\end{proof}
\end{proposition}

\par {\bf A).} Let $r=1$. There are two opportunities.
\par (I). Assume that $S'=\widetilde S$. Then $S=\PP(1,1,n)$ and $C=X_1$ is a
generator of cone. We obtain the cases $2-8(0)$, $2-9(0)$, $2-10(0)$.

\par (II). Assume that $S'\ne \widetilde S$ and $\widehat S=\FFF_n$. By proposition
\ref{vidc} we have $\widehat B_i\sim l_iE_0$ for all $i$. By the equality
$2=a_1+\sum l_ib_i$ and lemma \ref{vida2}
we get either $l_i=1$ for all $i$ and
$(b_1,b_2)=(\frac12,\frac23)$, $(\frac12,\frac34)$,
$(\frac12,\frac45)$, $(\frac12,\frac56)$, $(\frac23,\frac23)$,
$(\frac23,\frac34)$, or $\widehat B_1\sim 2E_0$ and $b_1=2/3$.
Let $l_i=1$ for all $i$.
Then $n+2=n\sum b_i+b+\gamma$ and there are two variants to extract a required curve.

\par (1). Let $\big(\CC^2,(2-n(b_1+b_2-1)-b)\{x=0\}+b_i\{y=0\}\big)$ and $\vartheta=1$.
\\
Thus $\alpha\ge \beta+1$ and similarly to lemma \ref{weight} we have
$$
\alpha(n(b_1+b_2-1)+b-1)+\beta(1-b_i)=1.
$$

Then $(n,b_i,b_1,b_2,\alpha,\beta)=(2,\frac12,\frac12,\frac23,2,1)$,
$(2,\frac23,\frac12,\frac23,3,1)$,
$(2,\frac34,\frac12,\frac34,2,1)$.
We obtain the cases $6-11(0)$, $9-5(0)$, $6-12(0)$.

\par (2). Let $\big(\CC^2,(2-n(b_1+b_2-1)-b)\{x=0\}+b_1\{y+x^k\}+b_2\{y=0\}\big)$,
where $1\le k\le n$.
Thus $\alpha\ge \beta+1$ and similarly to lemma \ref{weight} we have
$$
\alpha(n(b_1+b_2-1)+b-1)-\beta(b_1+b_2-1)=1/\vartheta.
$$

Take a blow-up with weights (1,1). Then the requirement the pair to be
$\frac17$-log terminal implies
$b_1+b_2+1/7<n(b_1+b_2-1)+b$.
We obtain $\vartheta=1$ and $(n,b_1,b_2,\alpha,\beta)=(2,\frac12,\frac34,3,1)$,
$(2,\frac12,\frac34,4,3)$, $(2,\frac12,\frac34,5,4)$,
$(2,\frac12,\frac45,3,2)$, $(2,\frac12,\frac45,4,3)$,
$(2,\frac12,\frac56,3,2)$, $(2,\frac23,\frac23,3,2)$,
$(2,\frac23,\frac34,2,1)$, $(3,\frac12,\frac23,3,1)$,
$(3,\frac12,\frac23,3,2)$,
$(3,\frac12,\frac23,4,3)$, $(3,\frac12,\frac34,2,1)$,
$(4,\frac12,\frac23,2,1)$.
In all cases $k=1$, except the last one in which $k\le 2$.
We obtain the cases $9-6(0)$, $11-3(0)$, $16(0)$, $8-6(0)$,
$11-3(0)$, $8-6(0)$, $8-7(0)$, $6-13(0)$, $14-3(0)$,
$13-4(0)$, $21-2(0)$, $7-5(0)$, $10-3(0)$.

\par  The case $\widehat B_1\sim 2E_0$ and $b_1=2/3$ was included in the case
from (2) in which $\widehat B_1$ is decomposed into two sections.

\par {\bf B).} Let $r=2$. There are two variants.
\par (I). Assume that $S'=\widetilde S$ and $\widehat S=\FFF_n$.
Then $\rho(S')=3$. Put $\widehat B=\sum b_i\widehat B_i$. By lemma
\ref{vida} $a_1+a_2\ge 6/7$. Hence $\widehat B\cdot f=2-a_1-a_2=\frac{k-1}k$,
where $k=2,3,4,5,6,\{\infty\}$.
From a structure of $\FFF_n$ we have $n+2=b+a_2n+\widehat B\cdot E_0$.
To be definite, assume that $\widetilde E_1^2\ge \widetilde E_2^2$
if $\Sing \widetilde S\cap \widetilde E_2=\emptyset$.
Taking into account proposition \ref{vidc} we have the following possibilities.

\par (1). $\widehat D=b\widehat C+a_1\widehat E_1+a_2\widehat E_2+\frac12\widehat B_1$,
where $n=2$ and $\widehat B_1\sim 2E_0$. It is possible that
$\widehat B_1$ is decomposed into
$\widehat B_1'$ and $\widehat B_1''$, where $\widehat B_1'\sim\widehat B_1''\sim E_0$.
Thus $a_2=1-b/2$. Considering all variants we obtain only one:

$$
\big(\CC^2,(1-b/2)\{x=0\}+\frac12\{(x+y^3)y=0\}\big)\ \text{and}\
(\vartheta,\alpha,\beta)=(1,7,2).
$$

In this variant $\widehat B_1$ is an irreducible divisor. It is the case
$17-2(0)$.

\par (2). $\widehat D=b\widehat C+a_1\widehat E_1+a_2\widehat E_2+
\frac{k-1}k\widehat B_1$, where $\widehat B_1\sim E_0+f$ and $k=2,3,4$.
Then $a_2=\frac1n(1+\frac1k-b)+\frac1k$.
\par
Let $k=2$. By proposition \ref{vidc} it follows that $n\le 4$ and $\widehat B_1$
intersects
$\widehat E_2$ at most at two different points.
\par
Assume that $n=2$ and $\widehat B_1$ is tangent $\widehat E_2$ with a multiplicity 3.
We have the pair $(\CC^2,(\frac54-\frac{b}2)\{x=0\}+\frac12\{x+y^3=0\})$.
The condition $\delta(S,D)=1$ implies $3b/2-9/4>-6/7$, i.e. $b>13/14$.
Similarly to lemma \ref{weight} we obtain $\vartheta=1$ and $(\alpha,\beta)=(9,2)$,
(11,3), (13,4), (16,5). These are the cases
$28(0)$, $31(0)$, $37(0)$, $49(0)$.

\par
Assume that $n=2$ and $\widehat B_1$ is tangent $\widehat E_2$ with a multiplicity 2.
Similarly we obtain $(\alpha,\beta)=(5,1)$ and $\vartheta=1$. It is the case
$12-2(0)$.

\par
Assume that $n=3$ and $\widehat B_1$ is tangent $\widehat E_2$ with a multiplicity 4.
Similarly the condition $\delta(S,D)=1$ implies $b>6/7$.
We obtain $\vartheta=1$ and $(\alpha,\beta)=(13,3)$, (17,4). These are the cases
$38(0)$, $41(0)$.

\par
Assume that $n=3$ and $\widehat B_1$ is tangent $\widehat E_2$ with a multiplicity 3.
There are no cases.

\par
Assume that $n=4$ and $\widehat B_1$ is tangent $\widehat E_2$ with a multiplicity 5.
Similarly we obtain $\vartheta=1$ and (11,2), (16,3). These are the cases
$32(0)$ and $50(0)$.
The remaining variants of intersection $\widehat B_1$ with $\widehat E_2$
are not realized.
\par
Let $k=3$. By proposition \ref{vidc} it follows that $n\le 3$.
\par
Assume that $n=2$ and $\widehat B_1$ is tangent $\widehat E_2$ with a multiplicity 3.
Similarly we obtain $\vartheta=1$ and (7,2). It is the case $17-3(0)$.
\par
The remaining variants of intersection $\widehat B_1$ and $\widehat E_2$
are not realized.
The case $n=3$ is not realized too.

\par
Let $k=4$. By proposition \ref{vidc} it follows that $n=2$.
Similarly the curve $\widehat B_1$ must be tangent $\widehat E_2$
with a multiplicity 3.
We obtain $\vartheta=1$ and (4,1). It is the case $47-2(0)$.

\par (3). $\widehat D=b\widehat C+a_1\widehat E_1+a_2\widehat E_2+
\frac{k-1}k\widehat B_1+\frac{k-1}k\widehat B_2$, where
$\widehat B_1\sim E_0$ and $\widehat B_2\sim f$.
This case was included in the answer of point (2).

\par (4). $\widehat D=b\widehat C+a_1\widehat E_1+a_2\widehat E_2+
\frac{k-1}k\widehat B_1$, where
$\widehat B_1\sim E_0$ and $2\le k\le 6$.
Then $a_2=\frac1k+\frac1n(2-b)$. Let $\nu$ be an intersection multiplicity of
$\widehat B_1$ and $\widehat E_2$ at the extraction point. By proposition \ref{vidc}
it follows either $\nu=n$ or $n-1$. In the case $\nu=n$ a discrepancy is equal to
$b-2<-1$ after a blow-up with weights $(n,1)$. We get a contradiction.
Therefore $\nu=n-1$, $\vartheta=1$ and $\beta=1$.
Taking into account $\widetilde E_1^2>\widetilde E_2^2$ and $\delta(S,D)=1$ we obtain
$(n,k,\alpha,\beta)=(2,3,5,1)$,
(3,2,7,1). These are the cases $12-3(0)$, $23(0)$.

\par (5). $\widehat D=b\widehat C+a_1\widehat E_1+a_2\widehat E_2+
\frac{k-1}k\widehat B_1+\frac{l-1}l\widehat B_2$, where
$\widehat B_1\sim E_0$, $\widehat B_2\sim f$, $k\ne l$ and $l\ge 2$.
Then $a_2=\frac1k+\frac1n(1+\frac1l-b)$. There are two variants:
the extraction takes place either at the point lying on
$\widehat B_2$, or at the point not lying on
$\widehat B_2$. In the first variant  the condition
$\delta(S,D)=1$ and proposition \ref{vidc} imply $\vartheta=1$, $\beta=1$
and the intersection multiplicity of
$\widehat B_1$ and $\widehat E_2$
is equal to $n-1$  at the extraction point.
Proving by exhaustion this case is impossible.
Similarly, in the second case we have that
$\vartheta=1$, $\beta=1$  and the intersection multiplicity is equal to
$n$. This case is impossible too.

\par (II). Assume that $S'\ne \widetilde S$, $E_2'\cap\Sing S'\ne\emptyset$
and $\widehat S=\FFF_n$. Then $\rho(\widetilde S)=4$ and two curves are extracted.
To be definite, assume that $a_1\le a_2$.
By proposition \ref{vidc} it follows that $\widehat B_i\sim l_iE_0$ and the curves
are extracted with $\vartheta=1$.
By lemma \ref{vida2} we have $\sum b_i=2-a_1-a_2\le 6/7$. Thus
$\widehat D=b\widehat C+a_1\widehat E_1+a_2\widehat E_2+
\frac{k-1}k\widehat B_1+\gamma\widehat\Gamma,$ where
$\widehat B_1\sim E_0$, $k=2,3,4,5,6$ if $n=2$, and $k=2$ if $n=3$.
\par (1). Let $n=2$. Then $4=2a_2+\gamma+b+2(\frac{k-1}k)$. Hence by lemma
\ref{vida2} we obtain $a_2\le \frac37+\frac1k$.
Put $P=\widehat \Gamma\cap \widehat E_2$. There are two possibilities.
\par
Let $(\widehat B_1\cdot \widehat E_2)_P=1$. Another intersection point of
$\widehat B_1$ and $\widehat E_2$ is denoted by $Q$.
When taking a blow-up at the point
$Q$ the self-intersection index of $\widehat E_2$ must be decreased at least by
3. Therefore, taking a blow-up with weights (3,1) we get $k=2$. Thus,
the extraction at the point
$Q$ can take place by a blow-up with weights
$(m,1)$, where $m=3,4$
(if $m\ge 5$ then $\delta(S,D)=2$). Then $a_2=1-\frac1{2m}$. Hence
$\gamma=1+\frac1m-b$. Let the blow-up at the point
$P$ has the weights $(\alpha,\beta)$. Then
$$
b=\frac{1+\beta(\frac12-\frac1{2m})}{\alpha}+\frac1m.
$$

Hence $(m,\alpha,\beta)=(3,3,2)$, (4,2,1). These are the cases
$18-2(0)$, $18-3(0)$.

\par
Let $(\widehat B_1\cdot \widehat E_2)_P=2$.
Then taking a blow-up with weights (2,1) at the point
$P$ we get a discrepancy $-2+b<-1$, a contradiction.

\par (2). Let $n=3\ (k=2)$. Similarly $a_2\le \frac{11}{14}$.
Put $P=\widehat \Gamma\cap \widehat E_2$.
If $(\widehat B_1\cdot \widehat E_2)_P=3$ then we obtain the same contradiction
with a log canonicality.
If $(\widehat B_1\cdot \widehat E_2)_P=2$ then the intersection multiplicity
is equal to 1 at another point $Q$.
When taking the extraction at the point
$Q$ the self-intersection index of
$\widehat E_2$ must be decreased at least by
4. Then we have a contradiction with
$a_2\le \frac{11}{14}$.
Thus $(\widehat B_1\cdot \widehat E_2)_P=1$ and
$(\widehat B_1\cdot \widehat E_2)_Q=2$. Similarly, taking a blow-up with weights
$(m,1)$ at the point $Q$ (here $m\ge 4$) we get
$$
a_2=1-\frac1m,\
\gamma=\frac12-b+\frac3m,\
b=\frac{1+\beta(\frac12-\frac1m)}{\alpha}+\frac3m-\frac12.
$$

Hence $(m,\alpha,\beta)=(4,2,1)$.
It is the case $18-4(0)$.

\par {\bf C).} Let $r=3$. To be definite, assume $a_1\le a_2\le a_3$. Then
$2=a_1+a_2+a_3+\frac{k-1}k$. Hence $a_1\le \frac13+\frac1{3k}$.
From a structure of
$\FFF_n$ it follows that
$n+2\ge (2-a_1)n+b+\gamma\ge (\frac53-\frac1{3k})n+\frac67+\gamma$.
Hence $\frac87-\gamma\ge (\frac23-\frac1{3k})n$. Therefore $(n,k)=(2,1)$, (2,2),
(2,3), (3,1) if $\gamma=0$ and $(n,k)=(2,1)$
if $\gamma\ne 0$. Also in all cases
$\sum b_i\widehat B_i=b_1\widehat B_1\sim_{\QQ}\frac{k-1}kE_0$.
\par
Let $(n,k,\gamma)=(2,1,0)$. The condition $K_{\widehat S}+\widehat D$ to be log
canonical implies that
$\widehat E_2$ intersects $\widehat E_3$ at two different points
$P_1$ and $P_2$. By proposition \ref{vidc}
we must take the blow-ups with weights $(m_i,1)$ at these points, where
$m_i\ge 3$. Moreover
$\vartheta=1$. The condition $\delta(S,D)=1$ implies
$m_i\le 4$. Thus

\begin{equation}
\left\{
\begin{array}{l}
m_1-m_1a_2-a_3=0  \\
m_2-m_2a_3-a_2=0  \\
a_2+a_3=2-\frac{b}2,\ b\in [\frac67,1).  \\
\end{array}
\right.
\end{equation}

Hence $(m_1,m_2)=(3,4)$ and $b=\frac{10}{11}$. It is the case
$55(0)$.

\par Let $\gamma=0$ and $(n,k)=(2,2)$ or (2,3). Taking into account $\delta(S,D)=1$
we have two possibilities: either $\widehat E_2$ intersects
$\widehat E_3$ at two different points $P_1$, $P_2$ and
$\widehat B_1$ passes through $P_1$, $P_2$, or
$\widehat E_2$ intersects
$\widehat E_3$ at the point $P$ with a multiplicity 2 and
$(\widehat B_1\cdot\widehat E_i)_P=1$, where $i=2,3$.
Writing out the equation system as
(3) we get that the cases are not realized.

\par Let $(\gamma,n,k)=(0,3,1)$.
Similarly $\widehat E_2$ must intersect
$\widehat E_3$ at the first point with a multiplicity 2 and at another point
with a multiplicity 1. As before this case is not realized.

\par Let $\gamma>0$, $(n,k)=(2,1)$ and
$\Sing \widetilde S\cap \widetilde E_i\ne \emptyset$, where $i=2,3$.
Then $\widehat E_2$ intersects
$\widehat E_3$ at two different points $P_1$ and $P_2$, and
$\widehat \Gamma$ passes through $P_2$. Therefore one curve is extracted
at the point
$P_1$  and two curves are extracted at the point
$P_2$. Thus, at the point $P_2$ we have the pair
$(\CC^2,\gamma\{x=0\}+a_2\{y=0\}+a_3\{x+y=0\})$.
Write out the equation system as (3). Then the extractions are the blow-ups
with weights (2,1) and (1,3) at the point $P_2$,
and a blow-up with weights (3,1) at the point $P_1$.
It is the case $56(0)$.

\par
\section{\bf Case ${\bf p_a(C)=0}$ and ${\bf \widetilde C^2=-1}$.}

Let $p_a(C)=0$ and $\widetilde C^2=-1$.
Let $P_1$, \ldots, $P_r$ be the singularities of $S$, which lie on $C$.
Recall that they are the cyclic singularities
$\frac1{n_i}(q_i,1)$. Put $B=\sum b_iB_i$ and
$n_1\le n_2\le\ldots \le n_r$.
By proposition \ref{CC2} we have
\begin{equation}\label{eq1}
\frac{-2+\sum_{i=1}^r\frac{n_i-1}{n_i}+B\cdot C}{1-b}-\sum_{i=1}^r\frac{q_i}{n_i}=-1.
\end{equation}

\begin{proposition} $r=2,3$.
\begin{proof} Since $C^2>0$ then by (\ref{eq1}) it follows that
$r\ge 2$. Let $r\ge 4$.
Assume that $r=4$ and $n_1=n_2=n_3=n_4=2$. Then $B\ne 0$. Thus
$B\cdot C\ge 1/4$. A contradiction with (\ref{eq1}).
\par
Assume that $r\ge 5$ or $n_r\ge 3$.
Then by (\ref{eq1}) it immediately follows that $r=4$, $n_1=n_2=n_3=2$,
$n_4=3$, $q_4=2$,
$b=\frac67$ and $B=0$. Let $S'\to S$ be a partial resolution of $S$ along $C$ taking
at the points $P_1$, $P_2$ and $P_3$, and $S'\to \widetilde S$ be a contraction of
proper transform of $C$.
Then $\rho(\widetilde S)=3$, $\widetilde E_i\cdot \widetilde E_j=3$ if
$i\ne j$ and $\widetilde E_i^2=1$.
Thus $K_{\widetilde S}+\sum_{i=1}^3\frac37\widetilde E_i\equiv 0$.
A linear system $|\widetilde E_1|$ gives a birational morphism
$\widetilde S\to \PP^2$ by proposition \ref{free2}. The images of $\widetilde E_i$
are denoted by
$\widehat E_i$. Then $\widehat E_1\sim \OO_{\PP^2}(1)$ and
$\widehat E_2\sim \widehat E_3\sim \OO_{\PP^2}(3)$.
Since $\widehat E_2^2=\widehat E_3^2=9$ then taking the extraction
we must decrease the
self-intersection index by 8. It is impossible to make it.
For example, let us consider the following case:
$(\CC^2,\frac37\{xy=0\}+\frac37\{(x+y^2)(y+x^2)=0\})$.
Taking a blow-up with weights
$(4,1)$ we obtain a discrepancy
$\frac47>0$, a contradiction.
\end{proof}
\end{proposition}

\par {\bf A).} Let $r=2$. Then

\begin{equation}\label{eq2}
B\ne 0 \ \ \ \text{and} \ \ \
\frac{B\cdot C-\frac1{n_1}-\frac1{n_2}}{1-b}-\frac{q_1}{n_1}-\frac{q_2}{n_2}=-1.
\end{equation}

Since $C^2>0$ then a divisor $B$ must intersect a curve $C$ at the smooth point of
$S$ by proposition \ref{vidc} .
From (\ref{eq2}) it follows that $B\cdot C<1$.
In particular, $B$ intersects $C$ at the single smooth point of $S$.
Consider case by case all intersection variants of
$B$ and $C$.
\par (I). Let $B$ intersects $C$ at the single point. Put
$B=\frac{k-1}kB_1$. From (\ref{eq2}) we get the following possibilities.

\renewcommand{\arraystretch}{1.5}
\begin{center}
\begin{tabular}{|c|c|c|c|c|}
\hline
No. & $k$ & $(n_1,q_1)$ & $(n_2,q_2)$ &$b$\\
\hline
1 & 2 & (3,2) & (7,4), (7,5), (7,6), (8,5), (8,7),&
$\frac9{10}, \frac{15}{16}, \frac{21}{22},\frac67, \frac{12}{13},$\\
&&&(9,7), (9,8), (10,9), (11,10), (12,11)&
$\frac78, \frac9{10}, \frac{15}{17}, \frac{33}{38}, \frac67$\\
\hline
2 & 2 & (3,1) & (7,6)& $\frac78$\\
\hline
3 & 2 & (4,3) & (5,3), (5,4), (6,5)&
$\frac67, \frac{10}{11}, \frac67$\\
\hline
4 & 3 & (2,1) & (7,5), (7,6), (8,7), (9,8)&
$\frac89, \frac{14}{15}, \frac89,\frac67$\\
\hline
\end{tabular}
\end{center}
\vspace{1mm}
\par
Consider the first case from the table. Let $S'\to S$ be a partial resolution of
$S$ along $C$ taking at the point $P_2$ and $\psi\colon S''\to S'$ be a partial
resolution of $S'$ along $E_2'$ and $\Exc\psi=\Gamma''$.
Let $S''\to \widetilde S$ be a contraction of proper transform of
$C$. Then
$\widetilde B_1\cdot\widetilde E_2=3$, $\widetilde E_2^2=1$ (see the table).
A linear system
$|\widetilde E_2|$ gives a birational morphism $\widetilde S\to \widehat S=\PP^2$
by proposition \ref{free2}
and
$\widehat E_2\sim \widehat \Gamma\sim\OO_{\PP^2}(1)$, $\widehat B_1\sim\OO_{\PP^2}(3)$.
Since $\rho(\widetilde S)=2$ then only one curve can be extracted
at the single intersection point of
$\widehat\Gamma$ and $\widehat B_1$. Denote this point by
$P$. Thus
$K_{\PP^2}+\frac12\widehat B_1+\gamma\widehat\Gamma+a_2\widehat E_2\equiv 0$.
At the point of tangency of $\widehat E_2$ and $\widehat B_1$
it must be extracted the curve with a discrepancy
$(-b)$. Hence $a_2=\frac12+\frac{b}3$ and
$\gamma=1-\frac{b}3$.

\par (1). Assume that we have the pair
$(\CC^2,(1-\frac{b}3)\{x=0\}+\frac12\{x+y^3=0\})$
in the neighborhood of $P$. Similarly to lemma \ref{weight}
we obtain $(\alpha,\beta)=(5,1)$, (7,2). These are the cases $25-2(-1)$,
$33(-1)$.

\par (2). Assume that we have the pair
$(\CC^2,(1-\frac{b}3)\{x=0\}+\frac12\{x^2+y^3=0\})$
in the neighborhood of $P$. We obtain a contradiction with log
canonicality.

\par (3). Assume that we have the pair
$(\CC^2,(1-\frac{b}3)\{x=0\}+\frac12\{y(x+y^2)=0\})$
in the neighborhood of $P$.
The condition
$\delta=1$ gives $b>\frac{27}{28}$.
By the same argument as in lemma
\ref{weight} we obtain that the possibility is not realized.

\par
Consider the second case from the table. Let $S'\to S$ be a partial resolution
$S'$ along $C$ taking at the point $P_1$ and
$S'\to \widetilde S$ be a contraction of proper transform of $C$. Then
$\rho(\widetilde S)=1$ and $\widetilde E_1^2=4$. Since
$\widetilde B_1\cdot \widetilde E_1=7$
then we have $\widetilde S=\PP(1,1,4)$ and $\widetilde B_1=X_7$
by proposition \ref{free2}. At the cone vertex
$K_{\widetilde S}+\frac12\widetilde B_1$ is not (1/7)-log terminal divisor, a
contradiction.

\par
Consider the third case from the table. The surface $\widetilde S$ is constructed in
the same way as in the first case from the table.
We have $\widetilde S=\PP(1,1,2)$ and
$\widehat E_2=X_2$, $\widehat \Gamma=X_1$, $\widehat B_1=X_4$.
Also $\widehat \Gamma$ intersects $\widehat B_1$ at the single point
$P$. Thus
$K_{\PP}+\frac12\widehat B_1+(\frac12+\frac{b}4)\widehat E_2+(1-\frac{b}2)
\widehat\Gamma\equiv 0$.
\par
If $P$ is a not cone vertex then similarly to lemma
\ref{weight} we obtain that the possibility is not realized.

\par
Let $P$ be a cone vertex and $\overline{B}_1$ and $\overline{\Gamma}$ are
the proper transforms
of $\widehat B_1$ and $\widehat\Gamma$ on a minimal resolution of $\PP(1,1,2)$.
Then $\overline{B}_1\sim E_0+2f$. There are two variants.

\par (1). Let $\overline{B}_1\cap E_{\infty}\cap\overline{\Gamma}\ne \emptyset$
and $\overline{B}_1$ is tangent to $E_{\infty}$. A contradiction with $\delta=1$.

\par (2). Let $\overline{B}_1\cap E_{\infty}\cap\overline{\Gamma}=Q_1$
and $\overline{B}_1$ intersects $E_{\infty}$ at two different points.
At the point $Q_1$ we have the pair
$(\CC^2,(1-\frac{b}2)\{x=0\}+(1-\frac{b}4)\{y=0\}+\frac12\{x+y=0\})$.
The condition $\delta=1$ gives $b>6/7$. Hence $(\alpha,\beta)=(4,3)$. It is the case
$39(-1)$.

\par
Consider the fourth case from the table. Let $S'\to S$ be a partial resolution of
$S$ along $C$. Take a partial resolution $\psi_1\colon S''\to S'$ of
$S'$ along $E_2'$ and $E_2^{(2)}=\Exc \psi_1$.
Take a partial resolution $\psi_2\colon S'''\to S''$ of
$S''$ along $E_2^{(2)}$ and $\Gamma'''=\Exc \psi_2$.
Let $S'''\to\widetilde S$ be a contraction of proper transform of $C$ and be
a contraction of all exceptional curves over
$P_2$, except $\Gamma'''$.
Then
a linear system $|\widetilde E_1|$ gives a birational morphism
$\widetilde S\to \PP^2$ by proposition
\ref{free2}. In the notations as above
$K_{\PP^2}+\frac23\widehat B_1+\frac{b}2\widehat E_1+(1-\frac{b}2)\widehat \Gamma
\equiv 0$. Here $\widehat B_1\sim \OO_{\PP^2}(3)$,
$\widehat E_1\sim \widehat \Gamma\sim \OO_{\PP^2}(1)$ and
$\widehat\Gamma$ intersects $\widehat B_1$ with a multiplicity 2 at the point different
from
$\widehat B_1\cap\widehat E_1$.
Taking into account $\delta=1$ we obtain that the extraction can be
for the next pair only:
$(\CC^2, (1-\frac{b}2)\{x=0\}+\frac23\{x+y^2=0\})$ if $(\alpha,\beta)=(3,1)$.
It is the case $17-4(-1)$.

\par (II). Let $B$ intersects $C$ at two points.
To be definite, assume that the coefficient of irreducible divisor from
$B$, which intersects
$C$ at the smooth point of $S$ is equal to $\frac{k-1}k$, and
the coefficient of irreducible divisor from
$B$, which intersects
$C$ at the singular point of $S$ is equal to
$\frac{l-1}l$.
Now we refuse from the condition
$n_1\le n_2$. Then by (\ref{eq1}) and corollary 3.10 \cite{Sh1}
we have
\begin{equation}\label{eq3}
\frac{\frac{k-1}k-\frac1{n_1l}-\frac1{n_2}}{1-b}-\frac{q_1}{n_1}-\frac{q_2}{n_2}=-1.
\end{equation}

Hence $(k,l,n_1,q_1,n_2,q_2,b)=(2,2,4,3,3,2,\frac9{10})$, $(2,2,5,4,3,2,\frac67)$.
\\
Let $S'\to S$ be a partial resolution of $S$ along $C$ taking at the point
$P_1$ and $\psi\colon S''\to S'$ be a partial resolution of $S'$
along $E_1'$ and $\Exc \psi=\Gamma''$.
Consider a contraction
$S''\to \widetilde S$ of proper transform of $C$.
Then a linear system $|\widetilde E_1|$ gives a birational morphism
$\widetilde S\to \PP^2$.
Thus
$K_{\PP^2}+\frac12\widehat B_1+(\frac12+\frac{b}3)\widehat E_1+(1-\frac{b}3)
\widehat\Gamma\equiv 0$.
Here
$\widehat E_1\sim \widehat \Gamma\sim\OO_{\PP^2}(1)$, $\widehat B_1\sim\OO_{\PP^2}(3)$.
As in the first case of the table consider the pair
$(\CC^2,(1-\frac{b}3)\{x=0\}+\frac12\{x+y^3=0\})$.
Note that $\widetilde\Gamma\cdot\widetilde B_1\ne 0$ after the extraction.
Then
$(\vartheta,\alpha,\beta)=(1,5,2)$ and hence $B=\frac12B_1$.
It is the case $18-5(-1)$.

\par (III). Let $B$ intersects $C$ at three points.
Similarly to (\ref{eq3}) we have

\begin{equation*}
\frac{\frac{k-1}k-\frac1{n_1l_1}-\frac1{n_2l_2}}{1-b}-
\frac{q_1}{n_1}-\frac{q_2}{n_2}=-1.
\end{equation*}
There are no solutions.
\par {\bf B).} Let $r=3$. Immediately note that this variant is not
realized.
From (\ref{eq1}) it follows that $B$ can't intersect $C$ at the smooth point of $S$.
We suppose that the coefficients of irreducible divisors from
$B$, which pass through
the points $P_i$ are equal to $\frac{l_i-1}{l_i}$ respectively. Then

\begin{equation}\label{eq4}
\frac{1-\sum_{i=1}^3\frac1{n_il_i}}{1-b}-
\sum_{i=1}^3\frac{q_i}{n_i}=-1.
\end{equation}

To be definite, assume that $n_1l_1\le n_2l_2\le n_3l_3$.
From (\ref{eq4}) it follows that $l_1=1$ and
$(n_1l_1,n_2l_2,n_3l_3)=(2,3,m)$, where $m\ge 7$; $(2,4,m)$, where $5\le m\le 12$;
$(2,5,m)$, where $5\le m\le 7$; $(2,6,6)$; $(3,3,m)$, where $4\le m\le 6$;
$(3,4,4)$. Moreover, the next proposition directly follows from
(\ref{eq4}).

\begin{proposition} \label{third}
There is one of the following cases:
\begin{enumerate}
\item there exists $i\ge 2$ such that $n_i<2q_i;$
\item $(n_1,q_1,n_2,q_2,l_2)=(2,1,3,1,1);$
\item $(n_1,q_1,n_2,q_2,n_3,q_3,l_2,l_3,b)=(2,1,2,1,5,2,2,1,\frac78).$
\end{enumerate}
\end{proposition}

\begin{remark} All possibilities from (2),
except $(n_3,q_3,b)=(7,3,\frac{10}{11})$ satisfy
the condition of the case (1). All possibilities from (2) will be considered
together for the accuracy.
\end{remark}

\par (I). Consider the first case from proposition \ref{third}.
Depending on $n_1$
there are two variants.
\par (1). Assume that $n_1=2$. Let $S'\to S$ be a partial resolution of $S$ along
$C$ taking at the points $P_2$ and $P_3$, and $\psi_1\colon S''\to S'$
be a partial resolution of
$S'$ along $E_i'$ and $\Gamma''=\Exc \psi_1$.
The number $i$ from proposition \ref{third} is determined in the next table.
If $\Sing S''\cap \Gamma''\ne \emptyset$ then take a partial resolution
$\psi_2\colon S'''\to S''$ of
$S''$ along $\Gamma''$ and $\Upsilon'''=\Exc \psi_2$.
Take a contraction $S'''\to\widetilde S$ of proper transform of $C$.
For $m\gg 0$ a linear system $|m\widetilde E_i+\widetilde \Gamma|$ gives a
birational morphism
$\widetilde S\to \widehat S$, where $\widehat S=\FFF_n$.
We have $K_{\widehat S}+\widehat\Theta\equiv 0$, where
$\widehat\Theta=d_1\widehat E_i+d_2\widehat\Gamma+d_3\widehat \Upsilon+d_4\widehat E_j+
\widehat D$
and $j\ne i$, $j\ge 2$.
It is clear that $\widehat E_i\sim \widehat \Upsilon\sim f$,
$\widehat \Gamma=E_{\infty}$,
$\widehat E_j\sim 2E_0$. Except the possibility 6 from table
we always have $\widehat B=0$.
All possibilities are described in the next tables.

\begin{center}
\begin{tabular}{|c|c|c|c|c|c|c|c|}
\hline
No. & $i$ & $(n_i,q_i)$ & $\rho(\widetilde S/\widehat S)$
&$\widehat S$ & $(n_j,q_j)$ &$b$ & $l_2,l_3$ \\
\hline
1& 2& $(3,2)$& 1& $\FFF_2$& & & $l_2=1$\\
\hline
2& 2& $(4,3)$& 2& $\FFF_2$& & & $l_2=1$\\
\hline
3& 3& $(5,3)$& 1& $\FFF_3$& $(4,1)$ & $\frac67$ & $l_2=l_3=1$\\
\hline
4& 3& $(5,4)$, $(6,5)$ & 2& $\FFF_2$& $(4,1)$ & $\frac{10}{11}$, $\frac67$ &
$l_2=l_3=1$\\
\hline
5& 3& $(5,3)$& 1& $\FFF_3$& $(2,1)$ & $\frac{11}{12}$ & $l_2=2$, $l_3=1$\\
\hline
6& 3& $(3,2)$& 1& $\FFF_2$& $(2,1)$ & $\frac78$ & $l_2=l_3=2$\\
\hline
7& 3& $(n_i,n_i-1)$& 2& $\FFF_2$& $(2,1)$ &  & $l_2=2$, $l_3=1$\\
&&$5\le n_i\le 8$&&&&&\\
\hline
8& 2& $(5,3)$& 1& $\FFF_3$&  &  & $l_2=l_3=1$\\
\hline
9& & $(5,4)$, $(6,5)$& 2& $\FFF_2$&  &  & $l_2=1$\\
\hline
\end{tabular}
\end{center}

\vspace{2mm}
Let us remark that the remaining information in this table is not necessary because
it is not involved in the further proof.
Now we write out the divisor
$\widehat\Theta$, which corresponds to the possibilities from this table.
\begin{center}
\begin{tabular}{|c|c|}
\hline
No. & Divisor $\widehat\Theta$ $(K_{\widehat S}+\widehat\Theta\equiv 0)$\\
\hline
1 & $\frac{2b}3\widehat E_i+\frac{b}3\widehat\Gamma+(1-\frac{b}6)\widehat E_j$\\
\hline
2 & $\frac{3b}4\widehat E_i+\frac{b}2\widehat\Gamma+\frac{b}4\widehat\Upsilon+
(1-\frac{b}4)\widehat E_j$\\
\hline
3 & $\frac57\widehat E_i+\frac47\widehat\Gamma+\frac57\widehat E_j$\\
\hline
4 & $\frac8{11}\widehat E_i+\frac6{11}\widehat\Gamma+\frac4{11}\widehat\Upsilon+
\frac8{11}\widehat E_j$, $\frac57\widehat E_i+\frac47\widehat\Gamma+
\frac37\widehat\Upsilon+\frac57\widehat E_j$\\
\hline
5 & $\frac34\widehat E_i+\frac7{12}\widehat\Gamma+\frac{17}{24}\widehat E_j$\\
\hline
6 & $\frac34\widehat E_i+\frac58\widehat\Gamma+\frac{11}{16}\widehat E_j+
\frac12\widehat B_k$,
where $\widehat B_k\sim f$\\
\hline
7 & $\frac34\widehat E_i+(\frac32-b)\widehat\Gamma+
(\frac94-2b)\widehat\Upsilon+(\frac14+\frac{b}2)\widehat E_j$\\
\hline
8 & $(\frac15+\frac{3b}5)\widehat E_i+(\frac25+\frac{b}5)\widehat\Gamma+
(\frac45-\frac{b}{10})\widehat E_j$\\
\hline
9 & $\frac{4b}5\widehat E_i+\frac{3b}5\widehat\Gamma+\frac{2b}5\widehat\Upsilon+
(1-\frac{3b}{10})\widehat E_j$ \\
\hline
\end{tabular}
\end{center}

\vspace{2mm}
Using the condition $\delta(S,D)=1$ and the proof of lemma \ref{weight}
the reader will have no difficulty in showing that the required extractions are
absent.

\par (2). Assume that $(n_1,q_1)=(3,2)$. Then $l_2=l_3=1$.
Let $S'\to S$ be a partial resolution of $S$ along
$C$ taking at the points $P_2$ and $P_3$, and $\psi\colon S''\to S'$
be a partial resolution of
$S'$ along $E_3'$ and $\Gamma''=\Exc \psi$.
Consider a contraction $S''\to\widetilde S$ of proper transform of $C$.
Then a linear system $|\widetilde E_3|$ gives a birational morphism
$\widetilde S\to \PP^2$ and
$K_{\PP^2}+x\widehat E_2+(1+\frac{b}3-x)\widehat E_3+
(2-2x-\frac{b}3)\widehat\Gamma\equiv 0$,
where $\widehat E_2\sim \OO_{\PP^2}(3)$,
$\widehat E_3\sim\widehat\Gamma\sim \OO_{\PP^2}(1)$.
Depending on the singularity at the point
$P_2$ we have $x=\frac{2b}3$,
$\frac13+\frac{b}3$ or
$\frac{3b}4$. The corresponding singularities are
$\AAA_2$, $\frac13(1,1)$, $\AAA_3$.
It can easily be checked that the extraction of two required curves is absent.

\par (II). Consider the second case from proposition \ref{third}.
Let $S'\to S$ be a partial resolution of $S$ along
$C$ taking at the point $P_3$, and $\psi\colon S''\to S'$ be a partial resolution of
$S'$ along $E_3'$ and $\Gamma''=\Exc \psi$.
Consider a contraction $S''\to\widetilde S$ of proper transform of $C$.
Then $p_a(\widetilde E_3)=1$ and from (\ref{eq4}) it follows that
$\widetilde E_3^2=3,4$. Similarly to proposition \ref{third} it is proved that
$|\widetilde E_3|$ is not very ample linear system, except the case
$S''=S'$. In any case it gives a birational morphism
$\widetilde S\to \widehat S$, where $\widehat S$ is Del Pezzo surface
with Du Val singularities. Since the curve
$\widetilde E_3$ has a cusp then $\pi_1(\widehat S\setminus\Sing\widehat S)=0$
by remark 1.12 \cite{GZ}. Therefore, by lemma 6 \cite{MZ} it follows that
$\widehat S=S(\DDD_5)$ or $S(\EEE_6)$. Thus
$K_{\widehat S}+(\frac23+\frac{b}6)\widehat E_3+(\frac13-\frac{b}6)\cdot
(K_{\widehat S})^2\ \widehat\Gamma\equiv 0$.
Here $\widehat\Gamma$ is $(-1)$ curve on a minimal resolution of
$\widehat S$ since $\widehat E_3\cdot\widehat \Gamma=1$.
In the case $S''=S'$ the component of $B$ which passes through $P_3$ will be
the curve $\widehat\Gamma$.

\par 1).$\widehat S=S(\DDD_5)$. Then
$K_{\widehat S}+(\frac43-\frac{2b}3\widehat)\Gamma$ is not log canonical divisor
(see point {\bf F} from \S 4), a contradiction.

\par 2).$\widehat S=S(\EEE_6)$. Then
$K_{\widehat S}+(1-\frac{b}2\widehat)\Gamma$ is not log canonical divisor
(see point {\bf I} from \S 4), a contradiction.

\par (III). Consider the third case from proposition \ref{third}.
Let $S'\to S$  be a partial resolution of $S$ along
$C$ taking at the points $P_2$ and $P_3$, and $\psi\colon S''\to S'$
be a partial resolution of
$S'$ along $E_3'$ and $\Gamma''=\Exc \psi$.
Consider a contraction $S''\to\widetilde S$ of $E_3''$
and proper transform of $C$.
As in point (II) a linear system
$|\widetilde E_2|$ gives a birational morphism
$\widetilde S\to \widehat S=S(\DDD_5)$.
Thus $K_{\widehat S}+\frac{11}{16}\widehat E_2+\frac12\widehat B_1+\frac38\widehat
\Gamma\equiv 0$, where
$\widehat B_1$ is $(-1)$ curve on a minimal resolution of
$\widehat S$, and $\widehat \Gamma\sim_{\QQ}\frac12\widehat E_2$.
Since $-1=(K_{\widehat S}+\widehat\Gamma)\cdot\widehat \Gamma=
-2+\Diff_{\widehat\Gamma}(0)$ then $\Diff_{\widehat\Gamma}(0)=1$. By the classification
of two-dimensional singularities
(for example, see \cite[theorem 2.1.3]{PrLect})
it follows that at the point
 $\DDD_5$ we have

\begin{gather*}
\Big(\widehat S,(1/2)\widehat B_1+(3/8)\widehat\Gamma\Big)\an \Big(x^2+y^2z-z^4=0
\subset (\CC^3,0),\\ (1/2)\{y=x-z^2=0\}+(3/8)\{x=z=0\}\Big).
\end{gather*}

Taking a blow-up with weights $(4,3,2)$ the discrepancy of exceptional divisor
is equal to
$-9/8<-1$. We obtain a contradiction.

\end{document}